\documentclass[twoside,11pt]{article}
\usepackage[english]{babel}
\usepackage{amssymb, theorem, psfrag, epsfig, amsmath, bm}
\usepackage[nohead,a4paper]{geometry}
\theoremheaderfont{\bfseries\upshape}
\def\proof#1{{\bf #1.}}
\def\endproof{\hfill$\Box$}

\input cyracc.def 

\newcommand{\bb}{\mathbb}
\def \R {{\bb R}}

\def \Z {{\bb Z}}

\def \pr{\bb P}
\def \mean{{\bb E}}

\def \prdown{\pr^{\downarrow}}
\def \meandown{\mean^{\downarrow}}
\def \Dcal{{\mathcal D}}
\def \C{{\mathcal C}}

\def \oX{{\overline X}}

\newcommand{\iz}{{\rm \rlap Y\kern 2.2pt Y}}

\newcommand{\calC}{{\cal C}}

\newcommand{\e}[1]{{\rm e}^{{#1}}}
\newcommand{\refs}[1]{(\ref{#1})}

\renewcommand{\theequation}%
{\arabic{equation}}

\newcounter{mylistcnt}
\renewcommand{\themylistcnt}{{\rm({\roman{mylistcnt}})}}

\newcounter{zad}

\newtheorem{Th}{Theorem}[section]

\newtheorem{Prop}{Proposition}[section]
\newtheorem{Lemma}{Lemma}[section]
\newtheorem{Ex}{Example}[section]
\newtheorem{Cor}{Corollary}[section]
\newtheorem{Defin}{Definition}[section]
\newtheorem{D}{}[section]
\newtheorem{Rem}{Remark}[section]

\begin{document}

\title{Quasi-product forms for L\'evy-driven fluid networks}
\author{K.~D\c{e}bicki$^{1}$, 
A.~B.~Dieker$^{2, \,3}$, 
T. Rolski$^{1}$}
\maketitle

\stepcounter{footnote}\footnotetext{Mathematical Institute,
University of Wroc\l aw, pl. Grunwaldzki 2/4, 50-384 Wroc\l aw,
Poland} 
\stepcounter{footnote}\footnotetext{CWI, P.O. Box 94079, 1090 GB Amsterdam, the Netherlands}
\stepcounter{footnote}\footnotetext{University of Twente, P.O. Box 217, 7500 AE Enschede, the Netherlands}
\begin{abstract}
We study stochastic tree fluid networks driven by a multidimensional
L\'evy process.
We are interested in (the joint distribution of) the steady-state
content in each of the buffers, the busy periods, and the idle periods.
To investigate these fluid networks, we relate the above three quantities
to fluctuations of the input L\'evy process by
solving a multidimensional Skorokhod reflection problem.
This leads to the analysis of the distribution of the componentwise maximums,
the corresponding epochs at which they are attained,
and the beginning of the first last-passage excursion.
Using the notion of splitting times, we are able to find their Laplace transforms.
It turns out that, if the components of the L\'evy process are
`ordered', the Laplace transform has a so-called quasi-product form.

The theory is illustrated by working out special
cases, such as tandem networks and priority queues.

\vskip 0.2cm \noindent {\em Keywords:}
$n$-dimensional L\'evy process, splitting time, quasi-product form,
multidimensional Skorokhod problem, L\'{e}vy-driven fluid network, tree fluid network,
buffer content, busy period, idle period.
\\
\vskip 0.1cm \noindent AMS 2000 Subject Classification: Primary: 60K25; Secondary: 90B05, 60G51.
\end{abstract}

\section{Introduction.}
Prompted by a series of papers by Kella and Whitt 
\cite{kella:paralleldependent1993,kella:stationarityfeedfw1997,kellawhitt:tandem1992,kellawhitt:martingales1992},
there has been a considerable interest in multidimensional
generalizations of the classical storage model with nondecreasing
L\'evy input and constant release rate
\cite[Ch.~4]{prabhu:storage1998}. In the resulting networks, often
called {\em stochastic fluid networks}, the input into the buffers
is governed by a multidimensional L\'evy process. Recently,
motivated by work of Harrison and Williams on diffusion
approximations~\cite{harrisonwilliams:open1987,harrisonwilliams:feedforward1992},
the presence of product forms has been investigated
\cite{kella:stabilitynonproduct1996,kella:nonproduct2000,konstantopoulos:networks2004,piera:productform2005}.
Recall that the stationary buffer-content vector has a product form
if it has independent components, meaning that the distribution of
this vector is a product of the marginal distributions.

The results in these papers show that, apart from trivial cases, the
stationary buffer-content vector of stochastic fluid networks {\em
never} has a product form. Despite this `negative' result, we show
that it may still be possible to express the joint distribution of
the buffer content in terms of the marginal distributions. This fact
is best visible in the Laplace domain. For certain tandem queues,
for instance, the Laplace transform is a product that cannot be
`separated'; we then say that the buffer-content vector has a {\em
quasi-product form}.

\medskip
In the literature on stochastic fluid networks, there has been a focus on
the stationary buffer-content vector $W$ or one of its components.
Here, we are also interested in the stationary distribution
of vector of ages of the busy periods $B$ and idle periods $I$.
The age of a busy (or idle) period
is the amount of time that the buffer content has been positive (or zero)
without being zero (positive).
Knowing these, it is also possible to find the distribution of the remaining
length of the busy (or idle) period and the total length of these periods.

We are interested in $W$, $B$, and $I$ for a class of L\'evy-driven fluid networks
with a tree structure,
which we therefore call {\em tree fluid networks}.
Our analysis of these networks relies on a detailed study of a
related multidimensional Skorokhod reflection problem
(see, e.g., Robert \cite{robert:stochnetworks2003}).
Using its explicit solution, we relate the triplet of vectors $(W,B,I)$
to the fluctuations of a multidimensional L\'evy process $X$.
We also prove that the stationary distribution
of the buffer-content vector is unique.

\medskip
Since our analysis of fluid tree networks is based on fluctuations
of the process $X$, this paper also contributes to
fluctuation theory for multidimensional L\'evy processes.
Supposing that each of the components of $X$ drifts to $-\infty$,
we write $\oX$ for the (vector of) componentwise maximums of $X$, $G$ for the
corresponding epochs at which they are attained, and $H$
for the beginning of the first last-passage excursion.
Under a certain independence assumption, if the components of $G$ are `ordered',
we express the Laplace transform of $(\oX,G)$ in terms of
the transforms of the marginals $(\oX_j,G_j)$.
Since $X_j$ is a real-valued L\'{e}vy process,
the Laplace transform of $(\oX_j,G_j)$ is known if
$X_j$ has one-sided jumps;
see for instance Bertoin~\cite[Thm.~VII.4]{bertoin:levy1996}.

We also examine the distribution of $H$ under the measure $\prdown_k$,
which is the law of $X$ given that the process $X_k$ stays nonpositive.
There exists a vast body of literature on (one-dimensional) L\'evy processes
conditioned to stay nonpositive (or nonnegative), see the recent paper by
Chaumont and Doney~\cite{chaumontdoney:conditioned2005} for references.
Under the measure $\prdown_k$, we also find the transform of $(\oX,G)$.
As a special case, we establish the Laplace transform
of the maximum of a L\'evy process conditioned to stay below a subordinator,
such as a (deterministic) positive-drift process.

\medskip
By exploiting the solution of the aforementioned Skorokhod problem,
the results that we obtain for the process $X$ can be cast immediately
into the fluid-network setting.
For instance, the knowledge of $(\overline X,G)$
allows us to derive the Laplace transform of the
stationary distribution of $(W, B)$ in a tandem network and a
priority system if there are only positive jumps, allowing Brownian
input at the `root' station.
That is, we characterize the {\em joint} law of the
buffer-content vector and the busy-period vector.
With the $\prdown_k$-distribution of $H$, we
establish the transform of the idle-period vector $I$ for a
special tandem network. Our formulas generalize all explicit
results for tandem fluid networks that are known to date (in the
Laplace domain), such as those obtained by
Kella~\cite{kella:paralleldependent1993} and more recently by
D\c ebicki, Mandjes and van Uitert~\cite{debicki:tandem2005}. Most
notably, quasi-products appear in our formulas, even for idle
periods.

\medskip
To derive our results, we make use of the notion of {\em splitting times}.
These essentially allow us to reduce the problem to the one-dimensional case.
For real-valued Markov processes, splitting times have been introduced
by Jacobsen~\cite{jacobsen:splitting1974}.
Splitting times decompose (`split') a sample path
of a Markov process into two independent pieces. A full
description of the process {\em before} and {\em after} the splitting time
can be given.
However, since the splitting time is not necessarily a stopping
time, the law of the second piece may differ from the original law of the Markov process.
We refer to Millar~\cite{millar:zeroone1977,millar:pathdecomposition1978} for further
details, and to Kersting and Memi\c{s}o\v{g}lu~\cite{kersting:pathdecompositions2004}
for a recent contribution.

The idea to use splitting times in the context of stochastic networks is novel.
The known results to date are obtained with
It\^o's formula~\cite{konstantopoulos:networks2004}, a closely related
martingale~\cite{kellawhitt:martingales1992}, or differential
equations~\cite{piera:productform2005}.
Intuitively, these approaches all exploit a certain harmonicity.
However, the results of Kyprianou and
Palmowski~\cite{kyprianoupalmowski:martingale2005}
already indicate that there is a relation between these approaches and splitting.
Splitting has the advantage that it is insightful and that proofs are short.
Moreover, it can also be used for studying  more complicated systems~\cite{diekermandjes:markovadditive2006}.

\medskip
This paper is essentially divided into two parts.
In the first part, consisting of Sections~\ref{sec:splittingtimes}--\ref{sec:distrH},
we analyze the fluctuations of an $n$-dimensional L\'evy processes $X$.
The notion of splitting times is formalized in Section~\ref{sec:splittingtimes}.
These splitting times are first used to study the distribution of $(\oX, G)$ in
Section~\ref{sec:distrXG}, and then to analyze the distribution of $H$
under $\prdown_k$ in Section~\ref{sec:distrH}.
The second part of this paper deals with fluid networks.
Section~\ref{sec:skorokhod} ties these networks to fluctuations of $X$,
so that the theory of the first part can be applied
in Section~\ref{sec:tandem}.
Finally, in Appendix~\ref{app:prelimcpdPsdrift}, we derive some
results for compound Poisson processes with negative drift.
They are used in Section~\ref{sec:distrH}.

\section{Splitting times.}
\label{sec:splittingtimes}

This paper relies on the application of {\em splitting times}
to a multidimensional L\'evy process.
After splitting times have been introduced,
we study splitting at the maximum (Section~\ref{sec:splitmax})
and splitting at a last-passage excursion (Section~\ref{sec:splitlastexit}).

Throughout, let $X=(X_1,\ldots,X_n)^{'}$ be
an $n$-dimensional L\'evy process,
that is, a c\`adl\`ag process with stationary, independent increments
such that $X(0)=0\in \R^n$.
Without loss of generality, as in Bertoin~\cite{bertoin:levy1996},
we work with the canonical measurable
space $(\Omega, \mathcal F)=(D([0,\infty),\R^d\cup\{\partial\}), \mathcal B)$,
where $\mathcal B$ is the Borel $\sigma$-field generated by
the Skorokhod topology, and $\partial$ is an isolated point that serves
as a cemetery state.
In particular, $X$ is the coordinate process.
Unless otherwise stated, `almost surely' refers to $\pr$.
All vectors are column vectors.

The following assumption is used extensively throughout this paper:

\begin{description}
\item[{\bf D}] $X_k(t)\to-\infty$ almost surely, for every $k$.
\end{description}

We emphasize that a dependence between components is allowed.
In the sequel, $\overline X_k(t)$ (or $\underline X_k(t)$)
is shorthand for $\sup_{s\leq t} X_k(s)$ (or $\inf_{s\leq t} X_k(s)$).
Due to {\bf D}, $\overline X_k:=\overline X_k(\infty)$
is well-defined and almost surely finite for every $k$.
Furthermore, we write $\overline X=(\overline X_1,\ldots,
\overline X_n)^{'}$.

The following two definitions are key to further analysis.
The second definition is closely related to the first,
but somewhat more care is needed on a technical level.
Intuitively, for the purposes of this paper,
there is no need to distinguish the two definitions.

\begin{Defin}
We say that a random time $T$ is a {\em splitting time} for $X$
under $\pr$ if the two processes $\{X(t): 0\leq t\leq T\}$ and
$\{X(T+t)-X(T):t\geq 0\}$ are independent under $\pr$. We say that
$T$ is a {\em splitting time from the left} for $X$ under $\pr$ if
the two processes $\{X(t): 0\leq t< T\}$ and $\{X(T+t)-X(T-):t\geq
0\}$ are independent under $\pr$.
\end{Defin}

Note that if $X$ is a L\'evy process under $\pr$
with respect to some filtration $\mathcal F$ which includes the natural filtration,
any $\mathcal F$-stopping time $\tau$
is a splitting time for $X$ under $\pr$.
In fact, the L\'evy assumption implies that
$\{X(\tau+t)-X(\tau):t\geq 0\}$ is not only independent of
$\{X(t): 0\leq t\leq \tau\}$, but that it also
has the same distribution as $\{X(t):t\geq 0\}$.

We need some notions related to the initial behavior of
$X$. For $k=1,\ldots,n$, set $\overline R_k=\inf\{t>0:
X_{k}(t)=\overline X_{k}(t)\}$. Since $\{\overline X_k(t)-X_k(t):t\ge0\}$
is a Markov process under $\pr$ with respect to the filtration generated
by $X$ (see Proposition~VI.1 of \cite{bertoin:levy1996}), the
Blumenthal zero-one law shows that either $\overline R_k>0$ almost
surely (0 is then called {\em irregular} for $\{\overline X_k(t)-X_k(t):t\ge0\}$) or
$\overline R_k=0$ almost surely (0 is then called {\em regular} for
$\{\overline X_k(t)-X_k(t):t\ge0\}$). We also set $\underline R_k=\inf\{t>0:
X_k(t)=\underline X_k(t)\}$, and define regularity of 0 for
$\{X_k(t)-\underline X_k(t):t\ge0\}$ similarly as for $\{\overline X_k(t)-X_k(t):t\ge0\}$. If
$\overline R_k=0$ almost surely, we introduce
\[
\overline S_k=\overline S_k^X:=\inf\{t>0: X_k(t)\neq \overline X_k(t)\}.
\]
Again, either $\overline S_k=0$ almost surely (0 is then called
an {\em instantaneous point} for $\{\overline X_k(t)-X_k(t):t\ge0\}$) or
$\overline S_k>0$ almost surely (0 is then called a {\em holding point}
for $\{\overline X_k(t)-X_k(t):t\ge0\}$). One defines $\underline S_k$,
instantaneous points, and holding
points for $\{X_k(t)-\underline X_k(t):t\ge0\}$ similarly if $\underline R_k=0$.

\subsection{Splitting at the maximum under $\pr$.}
\label{sec:splitmax}
Let $G_k=G_k^X:=\inf\{t\geq 0: X_k(t)=\overline X_k
{\rm\,\, or\,\,} X_k(t-)=\overline X_k\}$
be the (first) epoch that $X_k$ `attains' its maximum,
and write $G=(G_1,\ldots,G_n)^{'}$.
Observe that $G_k$ is well-defined and almost surely finite for every $k$
by {\bf D}.

\begin{Lemma}
\label{lem:splittingpath}
Consider a L\'evy process $X$ that satisfies {\bf D}.
\begin{enumerate}
\item[(i)]
If $\overline R_k>0$ $\pr$-almost surely or $X_k$ is a compound Poisson process,
then $G_k$ is a splitting time for $X$ under $\pr$.
\item[(ii)]
If $\overline R_k=0$ $\pr$-almost surely but $X_k$ is not a compound Poisson
process, then $G_k$ is a
splitting time from the left for $X$ under $\pr$.
\end{enumerate}
\end{Lemma}
\proof{Proof}
We use ideas of Lemma~VI.6 of Bertoin~\cite{bertoin:levy1996}, who
proves the one-dimensional case under exponential killing.

We start with the first case, in which the ascending ladder set
is discrete. Set $\tau_0=0$ and define the stopping times
$\tau_{n+1} = \inf\{t>\tau_n: \overline X_{k}(t)>\overline X_{k}(t-)\}$ for $n>0$.
Write $N=\sup\{n:\tau_n<\infty\}$. Note that {\bf D}
implies that $N<\infty$ almost surely.

Let $F$ and $K$ be bounded functionals, and apply the Markov property
to see that for $n\in\Z_+$,
\begin{eqnarray*}
\lefteqn{\mean \left[F(X(t), 0\leq t\leq G_{k})
K(X(G_{k}+t)-X(G_{k}), t\geq 0); N=n\right]} \\&=&
\mean \left[F(X(t), 0\leq t\leq\tau_n) 1_{\{N\geq n\}}
K(X(\tau_n+t)-X(\tau_n), t\geq 0)
1_{\{\sup_{t\geq \tau_n}X_{k}(t)= X_{k}(\tau_n)\}}\right] \\
&=& \mean \left[F(X(t), 0\leq t\leq \tau_n) 1_{\{N\geq n\}}\right]
\mean\left[K(X(\tau_n+t)-X(\tau_n), t\geq 0)
1_{\{\sup_{t\geq \tau_n}X_{k}(t)= X_{k}(\tau_n)\}}\right]\\
&=& \mean \left[F(X(t), 0\leq t\leq \tau_n) 1_{\{N\geq n\}}\right]
\mean\left[K(X(t), t\geq 0)1_{\{\sup_{t\geq 0}X_{k}(t)= 0\}}\right].
\end{eqnarray*}
Summing over $n$ shows that the processes
$\{X(t): 0\leq t\leq G_{k}\}$ and $\{X(G_{k}+t)-X(G_{k}): t\geq 0\}$
are independent.

The argument in the case $\overline R_k=0$ is more technical, but essentially the same.
The idea is to discretize the ladder height structure, for which we
use the local time $\overline \ell_k$ at zero of the process
$\{\overline X_k(t)-X_k(t):t\ge0\}$; see Bertoin~\cite[Ch.~IV]{bertoin:levy1996} for definitions.
Note that $\overline \ell_{k}(\infty)<\infty$ almost surely by Assumption {\bf D}.

Therefore, we fix some $\epsilon>0$, and denote the integer part of
$\epsilon^{-1}\overline \ell_{k}(\infty)$ by
$n=\lfloor\epsilon^{-1}\overline \ell_{k}(\infty)\rfloor$.
A variation of the argument for $\overline R_k>0$
(using the additivity of the local
time) shows that $\{X(t):0\leq t\leq \overline \ell_{k}^{\,-1}(n\epsilon)\}$ and
$\{X(\overline \ell_{k}^{\,-1}(n\epsilon)+t)- X(\overline \ell_{k}^{\,-1}(n\epsilon)):t\geq 0\}$
are independent. According to \cite[Prop.~IV.7(iii)]{bertoin:levy1996},
$\overline \ell_{k}^{\,-1}(n\epsilon)\uparrow G_{k}$ as $\epsilon\downarrow 0$,
which proves the lemma.
\endproof

\subsection{Splitting at a last-passage excursion under $\prdown_k$.}
\label{sec:splitlastexit}
Let $H_k=H_k^X:=\inf\{t\geq 0: \sup_{s\geq t} X_k(s)\neq X_k(t)\}$
be the beginning of the first last-passage excursion,
and write $H=(H_1,\ldots,H_n)'$.

In this subsection, we study the splitting properties of $H_k$
for some fixed $k=1,\ldots,n$.
We suppose that 0 is a holding point for $\{X_k(t)-\underline X_k(t):t\ge0\}$, i.e., that
$\underline R_k=0$ and $\underline S_k>0$ $\pr$-almost surely.
Under this condition, the event $\{\overline X_k=0\}$
has strictly positive probability.
Therefore, one can straightforwardly define
the conditional law $\prdown_k$ of $X$ given $\overline X_k=0$.

It is our aim to investigate splitting of $H_k$ under $\prdown_k$,
but we only have knowledge of $X$ under $\pr$.
As a first step, it is therefore useful to
give a {\em sample path construction} of the law $\prdown_k$
on the canonical measurable space $(\Omega,\mathcal F)$.
For this, we define a process $X^{k\downarrow}$ by
\begin{equation}
\label{eq:defXdownarrow}
X^{k\downarrow}(t)=\left\{ \begin{array}{cl}
X(t) & {\rm if} \:\: t\in \left[\underline R^{(j)}_k,\underline S^{(j)}_{k}\right);\\
X(\underline R^{(j)}_{k})-X((\underline R^{(j)}_{k}+\underline S^{(j)}_{k}-t)-) &
{\rm if} \:\: t\in \left[\underline S^{(j)}_{k},\underline R^{(j)}_{k}\right),
\end{array} \right.
\end{equation}
where $\underline R_k^{(0)}=0$, and for $j\geq 1$,
\[
\underline S^{(j)}_k :=
\inf\left\{t> \underline R^{(j-1)}_{k}: \underline X_k(t)\neq  X_k(t)\right\}, \quad
\underline R^{(j)}_k := \inf\left\{t>\underline
S^{(j)}_{k}:\underline X_k(t)= X_k(t)\right\}.
\]
In other words, $X^{k\downarrow}$ is constructed from
the coordinate process $X$ by `reverting' the excursions of
$\{X_k(t)-\underline X_k(t):t\ge0\}$.

We have the following interesting lemma, which is the key to all results
related to $\prdown_k$.
For the random-walk analogue, refer to
Doney~\cite{doney:tanaka2005}.

\begin{Lemma}
\label{lem:alternativeprdown}
Consider a L\'evy process $X$ that satisfies {\bf D}.
If $\underline R_k=0$ and $\underline S_k>0$ $\pr$-almost surely, then
$X^{k\downarrow}$ has law $\prdown_k$ under $\pr$.
\end{Lemma}
\proof{Proof}
Observe that $\overline R_k>0$, and that the post-maximum process
$\{X(G_k + t)- X(G_k): t\geq 0\}$ has distribution
$\prdown_k$ (a proof of this uses similar
arguments as in the proof of Lemma~\ref{lem:splittingpath};
see Millar~\cite{millar:zeroone1977,millar:pathdecomposition1978}
for more details).

Fix some $q>0$, and let $e_q$ be an exponentially distributed random
variable, independent of $X$ (obviously, one must then enlarge
the probability space). 
The first step is to construct the law of
$\{X(G^q_k + t)- X(G^q_k): 0\leq t< e_q-G^q_k\}$,
where $G_k^q:=\inf\{t< e_q: X_k(t)=\overline X_k(e_q) {\rm\,\, or\,\,}
X_k(t-)=\overline X_k(e_q)\}$.
By the time-reversibility of $X$ \cite[Lem.~II.2]{bertoin:levy1996},
it is equivalent to construct the
law of $\{X(F^q_k)- X((F^q_k-t)-): 0\leq t< F^q_k\}$, where
$F_k^q:=\sup\{t< e_q: X_k(t)=\underline X_k(e_q) {\rm\,\, or\,\,}
X_k(t-)=\underline X_k(e_q)\}$.

To do so, we use ideas from
Greenwood and Pitman~\cite{greenwoodpitman:splitting1980}.
Let $\underline \ell_k$ be the local time
of $\{X_k(t)-\underline X_k(t):t\ge0\}$ at zero
(since $\underline R_k=0$, $\overline S_k>0$, we refer to
Bertoin~\cite[Sec.~IV.5]{bertoin:levy1996} for its
construction).
Its right-continuous inverse is denoted by $\underline \ell_k^{-1}$.
The $X$-excursion at local time $s$, denoted by $X^s$,
is the c\`adl\`ag process defined by
\[
X^s(u):= X\left(\left(\underline \ell_k^{-1}(s-)+u\right)\wedge\underline \ell_k^{-1}(s)\right)-
X\left(\underline \ell_k^{-1}(s-)-\right), \quad u\geq 0.
\]
If $\underline \ell_k^{-1}(s-)=\underline \ell_k^{-1}(s)$, then we let $X^s$ be $\partial$,
the zero function that serves as a cemetery.
Since $\{X^s:s>0\}$ is a c\`adl\`ag-valued Poisson point process
as a result of {\bf D},
one can derive (e.g., with the arguments of Lemma II.2 and Lemma~VI.2
of \cite{bertoin:levy1996}) that the process
\[
W:=\left\{W(s)=\left(D(s), X^s\right):s>0\right\}
\]
is time-reversible, where $D(s):=X\left(\underline \ell_k^{-1}(s)\right)$.
After setting $\sigma_q:=\underline \ell_k^{-1}(e_q)$, it can be seen that
this implies that
$\{(D(s), X^s):0<s<\sigma_q\}$ and
$\{(D(\sigma_q-)-D((\sigma_q-s)-), X^{\sigma_q-s}): 0<s<\sigma_q\}$
have the same distribution.
In other words, one can construct the law of
$\{X(F^q_k)- X((F^q_k-t)-): 0\leq t< F^q_k\}$ from the law
of $\{X(t):0\leq t<F^q_k\}$ by `reverting' excursions
as in (\ref{eq:defXdownarrow}).

It remains to show that this construction is `consistent' in the
sense of Kolmogorov, so that one can let $q\to0$ to obtain the claim.
For this, note that the family $\{\sigma_q\}$
can be coupled with a single random variable through
$\sigma_q=\underline \ell_k^{-1}(e_1/q)$.
\endproof

\medskip
We now study the splitting properties of $H_k$
using the alternative construction of $\prdown_k$
given in Lemma~\ref{lem:alternativeprdown}.
Since $\underline S_k^{(1)}$ is a $\pr$-stopping time with respect to the
(completed) natural filtration of $X$, the Markov
property of $X$ under $\pr$ with respect to this filtration
\cite[Prop.~I.6]{bertoin:levy1996} immediately
yields the following analogue of Lemma~\ref{lem:splittingpath}.

\begin{Lemma}
\label{lem:splittingpathH}
Consider a L\'evy process $X$ that satisfies {\bf D}.
If $\underline R_k=0$ and $\underline S_k>0$ $\pr$-almost surely, then
$H_k$ is a splitting time for $X$ under $\prdown_k$.
Moreover, it has an exponential distribution under $\prdown_k$.
\end{Lemma}

We remark that the construction and analysis of $\prdown_k$
is the easiest under the assumption that $\underline R_k=0$ and
$\underline S_k>0$ $\pr$-almost surely, which is exactly what we need
in the remainder.
A vast body of literature is devoted to the case $n=1$, and
the measure $\prdown_1$ is then studied under the assumption that $\overline R_1=0$.
This is challenging from a theoretical point of view,
since the condition that the process stays negative has $\pr$-probability zero.
Therefore, much more technicalities are needed to treat this
case. We refer to Bertoin~\cite{bertoin:decomposition1991} and
Doney~\cite{doney:tanaka2005} for more details.
See also Chaumont and Doney~\cite{chaumontdoney:conditioned2005}.

\section{The $\pr$-distribution of $(\overline X,G)$.}
\label{sec:distrXG}

The aim of this section is to find the Laplace transform of the distribution
of $(\overline X,G)$, assuming some additional structure on
the process $X$. Thus, in the sequel we write
$X_k \prec X_j$ if there exists some $K_{kj} >0$ such
that $X_j- K_{kj} X_k$ is nondecreasing almost surely.

\begin{Lemma}\label{l.gi}
Suppose the L\'evy process $X$ satisfies {\bf D}. If $X_k\prec
X_j$, then $G_k\le  G_j$.
\end{Lemma}
\proof{Proof}
First note that $G_k,G_j<\infty$ as a consequence of {\bf D}.
To prove the claim, let us assume instead that $G_j<G_k$
while $\hat X(t):= X_j(t)-C X_k(t)$ is nondecreasing for some arbitrary $C>0$.
Suppose that $X_k(G_k)=\oX_k$ and $X_j(G_j)=\oX_j$; the argument can be repeated
if, for instance, $X_k(G_k-)=\oX_k$.
The assumption $G_j<G_k$ implies that
\[
0\leq \hat X(G_k)-\hat X(G_j) = X_j(G_k)- \oX_j - C\left[
\oX_k- X_k(G_j)\right]\leq 0,
\]
meaning that $\oX_k = X_k(G_j)$. This contradicts $G_j<G_k$ in view of
the definition of $G_k$.
\endproof

\medskip
The following proposition expresses the distribution of $(\overline
{X}, G)$ in terms of those of $(X(G_k), G_k)$ and $(X(G_{k}-),
G_{k})$.
We denote the scalar product of $x$ and $y$ in $\R^n$ by $\langle x,y\rangle$,
and we write `cpd Ps' for `compound Poisson'.
Throughout this paper,
the expression $\prod_j \alpha_j \times \prod_j\beta_j\times \gamma $ should be read
as $\left(\prod_j\alpha_j\right)\times\left(\prod_j\beta_j\right)\times \gamma$.

\begin{Prop}
\label{prop:splitting}
Suppose that $X$ is an $n$-dimensional
L\'evy process satisfying {\bf D} and that $X_1\prec X_2\prec
\ldots\prec X_n$. Then for any $\alpha, \beta\in \R_+^n$,
\begin{eqnarray*}
\mean \e{-\langle \alpha, G\rangle-\langle \beta, \oX
\rangle}
&=& \mathop{\prod_{j=1}^{n-1}}_{\overline R_j>0 \hbox{\rm \scriptsize \, or $X_j$ cpd Ps}} \frac{\mean
\e{-\left[\sum_{\ell=j}^n\alpha_\ell\right] G_{j}-\sum_{\ell=j}^n\beta_\ell X_{\ell} (G_{j})}}
{\mean \e{-\left[\sum_{\ell=j+1}^n\alpha_\ell\right] G_{j}-
\sum_{\ell=j+1}^n\beta_\ell X_{\ell}(G_{j})}}\\
&&\mbox{}\times \mathop{\prod_{j=1}^{n-1}}_{\overline R_j=0,
\hbox{\rm \scriptsize \, $X_j$ not cpd Ps}} \frac{\mean
\e{-\left[\sum_{\ell=j}^n\alpha_\ell\right]
G_{j}-\sum_{\ell=j}^n\beta_\ell X_{\ell} (G_{j}-)}}
{\mean \e{-\left[\sum_{\ell=j+1}^n\alpha_\ell\right] G_{j}
-\sum_{\ell=j+1}^n\beta_\ell X_{\ell}(G_{j}-)}} \times \mean \e{-\alpha_n G_{n}
-\beta_n \overline X_{n}}.
\end{eqnarray*}
\end{Prop}
\proof{Proof}
First observe that the assumptions imply that the
terms $X_{\ell}(G_{j})$ and $X_{\ell}(G_{j}-)$ in the formula are
nonnegative for $\ell\geq j$, which legitimates the use of the
Laplace transforms. Remark also that $\overline R_i=0$ for $i>j$
whenever $\overline R_j=0$, i.e., for some deterministic
$i_0$ we have $\overline R_i>0$ for $i\le i_0$ and $\overline
R_i=0$ for $i> i_0$.

Let us first suppose that $\overline R_j>0$ or that $X_j$ is a compound Poisson process.
We prove that for $j=1,\ldots, n-1$,
\[
\mean \e{-\sum_{\ell=j}^n \alpha_\ell G_{\ell}
-\sum_{\ell=j}^n \beta_\ell \oX_{\ell}} =
\frac{\mean \e{-\left[\sum_{\ell=j}^n\alpha_\ell\right] G_{j}
-\sum_{\ell=j}^n\beta_\ell X_{\ell}(G_{j})}} {\mean \e{-\left[
\sum_{\ell=j+1}^n \alpha_\ell\right] G_{j}-\sum_{\ell=j+1}^n\beta_\ell
X_{\ell}(G_{j})}} \mean
\e{- \sum_{\ell=j+1}^n\alpha_\ell G_{\ell}-\sum_{\ell=j+1}^n \beta_\ell \oX_{\ell}}.
\]

The key observations are that $\oX_{j}= X_{j}(G_{j})$
and that $G_{{\ell}}\ge G_{j}$ almost surely
for $\ell=j,\ldots, n$ by Lemma~\ref{l.gi}.
The fact that $G_j$ is a splitting time by Lemma~\ref{lem:splittingpath}(i) then yields
\begin{eqnarray}
\lefteqn{\mean \e{-\sum_{\ell=j}^n \alpha_\ell
G_{\ell}-\sum_{\ell=j}^n \beta_\ell \oX_{\ell}}} \nonumber\\&=& \mean \e{
-\left[\sum_{\ell=j}^n \alpha_\ell\right] G_{j}-\sum_{\ell=j}^n\beta_\ell
X_{\ell}(G_{j})} \e{-\sum_{\ell=j+1}^n \alpha_\ell \left[G_{\ell}-G_{j}\right]
-\sum_{\ell=j+1}^n \beta_\ell \left[\oX_{\ell}-X_{\ell}(G_{j})\right]}\nonumber\\
&=&\mean \e{-\left[\sum_{\ell=j}^n
\alpha_\ell\right] G_{j}-\sum_{\ell=j}^n \beta_\ell X_{\ell}(G_{j})}
\mean \e{- \sum_{\ell=j+1}^n \alpha_\ell
\left[G_{\ell}-G_{j}\right]-\sum_{\ell=j+1}^n \beta_\ell
\left[\oX_{\ell}-X_{\ell}(G_{j})\right]}.
\label{eq:splittingwithoutmarks}
\end{eqnarray}
The latter factor, which is rather complex to analyze directly,
can be computed upon choosing $\alpha_j=\beta_j=0$ in the above display.

Repeating this argument for the case $\overline R_j=0$ yields with
Lemma~\ref{lem:splittingpath}(i), provided that $X_j$ is
not a compound Poisson process,
\[
\mean \e{-\sum_{\ell=j}^n \alpha_\ell G_{\ell}-\sum_{\ell=j}^n \beta_\ell
\oX_{\ell}} = \frac{\mean \e{- \left[\sum_{\ell=j}^n \alpha_\ell\right] G_{j}
-\beta_j\oX_{j}-\sum_{\ell=j+1}^n \beta_\ell
X_{\ell}(G_{j}-)}} {\mean
\e{- \left[\sum_{\ell=j+1}^n \alpha_\ell\right]
G_{j}-\sum_{\ell=j+1}^n \beta_\ell X_{\ell}(G_{j}-)}}
\mean \e{-\sum_{\ell=j+1}^n \alpha_\ell G_{\ell}-\sum_{\ell=j+1}^n \beta_\ell \oX_{\ell}}.
\]
It is shown in the proof of Theorem~VI.5(i) of
\cite{bertoin:levy1996} that $\oX_{j} = X_{j}(G_{j}-)$ almost
surely, and this proves the claim.
\endproof

\medskip
In the rest of this section,
the following assumption is imposed.

\begin{description}
\item[{\bf G}] For $j=1,\ldots,n-1$, we have
\begin{equation}
\label{eq:recursionXUpsilon}
X_{j+1}(t)=K_{j+1}X_{j}(t)+\Upsilon_{j+1}(t),
\end{equation}
where $(\Upsilon_2,\ldots,\Upsilon_n)$ are mutually independent
nonnegative subordinators and $K_2,\ldots,K_n$ are strictly positive.
\end{description}

Note that Assumption {\bf G} implies $X_1\prec X_2\prec
\ldots\prec X_n$.
Moreover, it entails that for $j=1,\ldots,n-1$ and $\ell\geq j$, we
have
\[
X_{\ell}(t) = K_j^\ell X_{j}(t) + \sum_{i=j+1}^\ell K_i^\ell
\Upsilon_{i}(t),
\]
where we have set $K_j^\ell=\prod_{i=j+1}^\ell
K_{i}$ and $K_j^j=1$. In other words, $X_\ell$ can be written as the sum of
$X_{j}$ and $\ell-j$ independent processes, which are all mutually
independent and independent of $X_{j}$.

The following reformulation of \refs{eq:recursionXUpsilon}
in terms of matrices is useful in Section~\ref{sec:tandem}.
Let $K$ be the upper triangular matrix with element $(i,i+1)$
equal to $K_{i+1}$ for $i=1,\ldots,n-1$, and zero elsewhere.
Also write $\Upsilon(t):=(\Upsilon_1(t),\ldots,
\Upsilon_n(t))'$, where $\Upsilon_1(t)=X_1(t)$.
Equation (\ref{eq:recursionXUpsilon}) is then nothing else than the
identity $X(t)=(I-K')^{-1}\Upsilon(t)$. The matrix $(I-K')^{-1}$ is
lower triangular, and
element $(i,j)$ equals $K_j^i$ for $j\geq i$.

The {\em cumulant} of the subordinator $\Upsilon_j(t)$
is defined as
\[
\theta^\Upsilon_j (\beta) := -\log \mean \e{-\beta \Upsilon_{j}(1)}
\]
for $\beta\ge0$ and $j=2,\ldots,n$.

The following theorem expresses the
joint Laplace transform of $(\oX,G)$ in terms of its marginal
distributions and the cumulants $\theta^\Upsilon$. 
However, except for trivial cases, the Laplace transform
is {\em not} the product of marginal Laplace transforms.
Still, it can be expressed in terms of these marginal transforms
in a product-type manner. We call this a {\em quasi-product} form.

\begin{Th}
\label{thm:resultindep}
Suppose that $X$ is an $n$-dimensional L\'evy process satisfying
{\bf D} and {\bf G}. Then for any $\alpha,\beta\in\R_+^n$,
the transform
$\mean \e{-\langle \alpha, G\rangle-\langle \beta, \oX\rangle}$
equals
\[
\prod_{j=1}^{n-1} \frac{\mean \e{-\left[\sum_{\ell=j}^n \alpha_\ell+\sum_{\ell=j+1}^n
\theta^\Upsilon_{\ell} \left(\sum_{k=\ell}^n K_\ell^k
\beta_k\right)\right] G_{j}-\left[\sum_{\ell=j}^n K_j^\ell
\beta_\ell\right] \overline X_{j} }} {\mean
\e{-\left[\sum_{\ell=j+1}^n \alpha_\ell+\sum_{\ell=j+1}^n \theta^\Upsilon_{\ell} \left(\sum_{k=\ell}^n
K_\ell^k \beta_k\right) \right] G_{j}
-\left[\sum_{\ell=j+1}^n K_j^\ell \beta_\ell\right] \overline X_{j}}}
\times \mean \e{-\alpha_n G_{n}-\beta_n \overline X_{n}}.
\]
\end{Th}
\proof{Proof}
Let $j$ be such that $\overline R_j>0$ or $X_j$ is compound Poisson.
By Assumption {\bf G}, we then have for $a\in\R_+$,
\begin{eqnarray*}
\mean \e{-a G_{j}-\sum_{\ell=j}^n \beta_\ell X_{\ell}(G_{j}) } &=&
\mean \e{-a G_{j}-\left[\sum_{\ell=j}^n K_j^\ell
\beta_\ell\right] X_{j}(G_{j}) -
\sum_{\ell=j+1}^n \left[\sum_{k=\ell}^n K_\ell^k \beta_k\right]
\Upsilon_{\ell}(G_{j})}\\&=& \mean
\left(\e{-a G_{j}-\left[\sum_{\ell=j}^n K_j^\ell \beta_\ell\right]
X_{j}(G_{j})
}\mean \left[\left.\e{-\sum_{\ell=j+1}^n
\left[\sum_{k=\ell}^n K_\ell^k \beta_k\right]
\Upsilon_{\ell}(G_{j})}\right| G_{j}\right]\right)\\
&=& \mean \e{-\left[a+\sum_{\ell=j+1}^n \theta^\Upsilon_{\ell}
\left(\sum_{k=\ell}^n K_\ell^k \beta_k\right)\right]
G_{j}-\left[\sum_{\ell=j}^n K_j^\ell \beta_\ell\right]
X_{j}(G_{j})}.
\end{eqnarray*}
The claim now follows from Proposition~\ref{prop:splitting} and the
fact that $X_{j}(G_{j})=\overline X_j$ almost surely.

If $\overline R_j=0$ but not a compound Poisson process,
the same argument gives the joint transform of
$\{X_{\ell}(G_{j}-):\ell=j,\ldots,n\}$ and $G_{j}$. In the resulting
formula, $X_{j}(G_{j}-)$ can be replaced by $X_{j}(G_{j})$ as
outlined in the proof of Theorem~VI.5(i) in Bertoin~\cite{bertoin:levy1996}.
\endproof

\medskip
The following corollary shows that
Theorem~\ref{thm:resultindep} not only completely characterizes the
law of $(\oX,G)$ under $\pr$, but also its law conditioned on
one component to stay nonpositive.
Indeed, let $\prdown_k$ be the law of
$\{X(G_k+t)-X(G_k): t\geq 0\}$ for $k=1,\ldots,n$; it can be checked that this
measure equals $\prdown_k$ as defined in Section~\ref{sec:splitlastexit}
in case $\underline R_k=0$ and $\underline S_k>0$ $\pr$-almost surely.
Note that $\prdown_k$ can be regarded as the law of $X$ given that $X_k$ stays nonpositive.

\begin{Cor}
For $\alpha,\beta\in\R_+^n$, we have
\[
\meandown_k \e{-\langle \alpha,\oX\rangle -\langle \beta, G\rangle}=
\prod_{j=k}^{n-1}
\frac{\mean \e{-\left[\sum_{\ell=j+1}^n \alpha_\ell+\sum_{\ell=j+2}^n
\theta^\Upsilon_{\ell} \left(\sum_{i=\ell}^n K_\ell^i
\beta_i\right)\right] G_{j+1}-\left[\sum_{\ell=j+1}^n K_{j+1}^\ell
\beta_\ell\right] \overline X_{j+1} }}
{\mean
\e{-\left[\sum_{\ell=j+1}^n \alpha_\ell+\sum_{\ell=j+1}^n \theta^\Upsilon_{\ell} \left(\sum_{i=\ell}^n
K_\ell^i \beta_i\right) \right] G_{j}
-\left[\sum_{\ell=j+1}^n K_j^\ell \beta_\ell\right] \overline X_{j}}}.
\]
\end{Cor}
\proof{Proof}
Directly from Theorem~\ref{thm:resultindep} and (\ref{eq:splittingwithoutmarks}).
\endproof

\medskip
In particular, this corollary characterizes the law of the maximum of a
L\'evy process given that it stays below a subordinator.
It provides further motivation for studying the law of the
vector $H$ under $\prdown_k$.


\section{The $\prdown_k$-distribution of $H$.}
\label{sec:distrH}
The aim of this section is to find the Laplace transform of the distribution
of $H$ under $\prdown_k$ under the assumption that 0 is a holding
point for $\{X_k(t)-\underline X_k(t):t\ge0\}$ under $\pr$.

We try to follow the same train of thoughts that led
us to the results in Section~\ref{sec:distrXG}.
This analogy leads to Proposition~\ref{prop:analog}, which does not
yet give the Laplace transform of the distribution
of $H$ under $\prdown_k$.
Therefore, we need an auxiliary result, formulated
as Lemma~\ref{lem:recursionidle},
which relies on Appendix~\ref{app:prelimcpdPsdrift}.
Finally, Proposition~\ref{prop:distrho} enables us to
find the Laplace transform of the distribution
of $H$ under $\prdown_k$.

As in the previous section, additional assumptions are imposed
on the L\'evy process $X$.
Here, they are significantly more restrictive.
The following Assumption {\bf H} plays a similar role in the
present section as Assumption {\bf G} in Section~\ref{sec:distrXG}.
Note that it implies $X_1\prec X_2\prec \ldots\prec X_n$.

\begin{description}
\item[{\bf H}]
Let $\Pi=\{\Pi(t):t\geq 0\}$ be a compound Poisson process with positive jumps only.
For each $j=1,\ldots,n$, we have \[
X_{j}(t)= 
\Pi(t)-c_j t,\]
where $c_j$ decreases strictly in $j$.
\end{description}

In the remainder of this section, we write $\lambda\in(0,\infty)$ for the intensity of jumps
of $\Pi$.
We also set $\rho^{(n)}_k:=\sup\{\underline R^{(j)}_k:\underline R^{(j)}_k\leq
\underline R_n^{(1)}\}$ and
$\sigma^{(n)}_k:=\sup\{\underline S^{(j)}_k:\underline S^{(j)}_k\leq
\underline R_n^{(1)}\}$. In particular,
$\rho_n^{(n)}=\underline R_n^{(1)}$ and $\sigma_n^{(n)}=\underline S_n^{(1)}$.
Also, we write for $\beta\geq 0$ and $i=1,\ldots,n$,
\[
\psi_i(\beta):=\log \mean \e{-\beta X_i(1)}
\]
for the {\em Laplace exponent} of $-X_i$.
Since we assume {\bf D}, $\psi_i$ is strictly
increasing on $\R_+$, see the proof of Corollary~VII.2 of
Bertoin~\cite{bertoin:levy1996}.
Therefore, we can define $\Phi_i$ as the
inverse of $\psi_i$. The function $\Phi_i$ plays an
important role in this section.

Recall that we used $n$ splitting times
to arrive at Proposition~\ref{prop:splitting}. Here, we only know that $H_k$
is a splitting time for $X$ under $\prdown_k$ (see Lemma~\ref{lem:splittingpathH}).
In general, however, $H_i$ ($i<k$) is not a splitting time under $\prdown_k$,
and the similarity with Proposition~\ref{prop:splitting} is lost.

\begin{Prop}
\label{prop:analog}
Suppose the L\'evy process $X$ satisfies {\bf D}.
For $\gamma\in\R_+^k$, we have
\begin{eqnarray*}
\meandown_k \e{-\sum_{j=1}^k \gamma_j H_j}&=&
\frac{\lambda}{\lambda+\sum_{j=1}^k \gamma_j}
\mean \e{-\sum_{j=1}^{k-1} \gamma_j \left(\rho_k^{(k)}-
\rho^{(k)}_j\right)}.
\end{eqnarray*}
\end{Prop}
\proof{Proof}
Lemma~\ref{lem:splittingpathH} yields
\begin{eqnarray*}
\meandown_k \e{-\sum_{j=1}^k \gamma_j H_j}&=&
\meandown_k \e{-\left(\sum_{j=1}^k \gamma_j\right) H_k}
\meandown_k \e{-\sum_{j=1}^{k-1} \gamma_j (H_j-H_k)}.
\end{eqnarray*}
In the discussion following (\ref{eq:defXdownarrow}),
we have seen that there is a simple sample-path correspondence between
the laws $\prdown_k$ and $\pr$. This yields immediately that $H_k$ is exponentially
distributed under $\prdown_k$ with parameter $\lambda$.
It also gives that the $\prdown_k$-distribution of
$\{H_j-H_k:j=1,\ldots,k-1\}$ is the same
as the $\pr$-distribution $\{\rho_k^{(k)}-\rho^{(k)}_j:j=1,\ldots,k-1\}$.
\endproof

\medskip
Motivated by the preceding proposition, we now focus on the calculation of the
distribution of the $\rho_k^{(k)}-\rho^{(k)}_j$
(that is, their joint Laplace transform).
For this, we apply results from Appendix~\ref{app:prelimcpdPsdrift}.

The following lemma is of crucial importance, as it provides a
recursion for the transform of
$\{\rho^{(i)}_{j+1}-\rho^{(i)}_j:j=1,\ldots,i-1\}$ and
$\{\rho^{(i)}_j-\sigma^{(i)}_j:j=1,\ldots, i\}$
in terms of the transform of the same family with superscript $(i-1)$.
The transforms of the marginals $\rho^{(i)}_i-\sigma^{(i)}_i$ and
$\rho^{(i-1)}_{i-1}-\sigma^{(i-1)}_{i-1}$ also appear in the
expression, but these transforms are known:
for $\gamma\geq 0$, $i=1,\ldots,n$ (cf.~the proof of Proposition~\ref{prop:excdistr}),
\begin{equation}
\label{eq:transformlengths}
\lambda \mean \e{-\gamma \left(\rho^{(i)}_i-\sigma^{(i)}_i\right)} =
\lambda+\gamma - c_i\Phi_i(\gamma).
\end{equation}

\begin{Lemma}
\label{lem:recursionidle}
Suppose that $X$ is an $n$-dimensional L\'evy process
satisfying {\bf D} and {\bf H}. Then for any $i=2,\ldots,n$,
$\beta\in\R_+^{i-1},\gamma \in \R_+^i$,
we have the following recursion:
\begin{eqnarray*}
\lefteqn{\mean \e{-\sum_{j=1}^{i-1} \beta_j \left(\rho^{(i)}_{j+1}-\rho^{(i)}_j\right)
-\sum_{j=1}^i \gamma_j \left(\rho^{(i)}_j-\sigma^{(i)}_j\right)}}
\\ &=&\frac{\beta_{i-1}+\lambda \mean \e{-\gamma_i
\left(\rho^{(i)}_i-\sigma^{(i)}_i\right)}}
{\beta_{i-1}+\lambda \mean \e{-\left[\left(\frac{c_{i-1}}{c_i}-1\right)
(\lambda+\beta_{i-1})+ \frac{c_{i-1}}{c_i}\gamma_i\right]
\left(\rho^{(i-1)}_{i-1}-\sigma^{(i-1)}_{i-1}\right)}}
\\&&\mbox{}\times
\mean \e{-\sum_{j=1}^{i-2} \beta_j \left(\rho^{(i-1)}_{j+1}-\rho^{(i-1)}_j\right)
-\sum_{j=1}^{i-2} \gamma_j \left(\rho^{(i-1)}_j-\sigma^{(i-1)}_j\right)
-\left[\left(\frac{c_{i-1}}{c_i}-1\right)
(\lambda+\beta_{i-1})+ \frac{c_{i-1}}{c_i}\gamma_i+\gamma_{i-1}\right]
\left(\rho^{(i-1)}_{i-1}-\sigma^{(i-1)}_{i-1}\right)}.
\end{eqnarray*}
\end{Lemma}
\proof{Proof}
Fix some $i=2,\ldots,n$, and consider the process
$X_{i-1}$ between $\sigma^{(i)}_i$ and $\rho_i^{(i)}$.
There are several excursions (at least one)
of the process $\{X_{i-1}(t)-\underline X_{i-1}(t):t\ge0\}$ away from 0 between
$\sigma^{(i)}_i$ and $\rho_i^{(i)}$, and
we call these excursions the $(i-1)$-subexcursions.
Each $(i-1)$-subexcursion contains excursions of the
processes $\{X_\ell(t)-\underline X_\ell(t):t\ge0\}$ for $\ell<i-1$; we call these the
$\ell$-subexcursions.
To each $(i-1)$-subexcursion, we assign $2i-4$ marks, namely two for
each of the $i-2$ types of further subexcursions. The first mark corresponds to
the length of the last $\ell$-subexcursion in the $(i-1)$-subexcursion, and
the second to the difference
between the end of the last $\ell$-subexcursion and the
end of the $(\ell+1)$-subexcursion.
Observe that these marks are
independent for every $(i-1)$-subexcursion between $\sigma^{(i)}_i$ and $\rho_i^{(i)}$,
and that their distributions
are equal to those of $\{\rho^{(i-1)}_\ell-\sigma^{(i-1)}_\ell: \ell=1,\ldots, i-2\}$
(the first marks)
and $\{\rho^{(i-1)}_{\ell+1}-\rho^{(i-1)}_\ell: \ell=1,\ldots, i-2\}$ (the second
marks).

\begin{figure}\centering
\psfrag{Xi-1}[Bl][Bl]{\huge $X_{i-1}$}
\psfrag{Z}[Bl][Bl]{\huge $Z$}
\psfrag{0}[Bl][Bl]{\huge $0$}
\psfrag{slope}[Br][Br]{\huge slope: $c_i-c_{i-1}$}
\psfrag{tau_-}[Br][Br]{\huge $\tau_-$}
\psfrag{rii}[Br][Br]{\huge $\rho_i^{(i)}$}
\psfrag{sii}[Br][Br]{\huge $\sigma_i^{(i)}$}
\resizebox{130mm}{!}{\includegraphics*{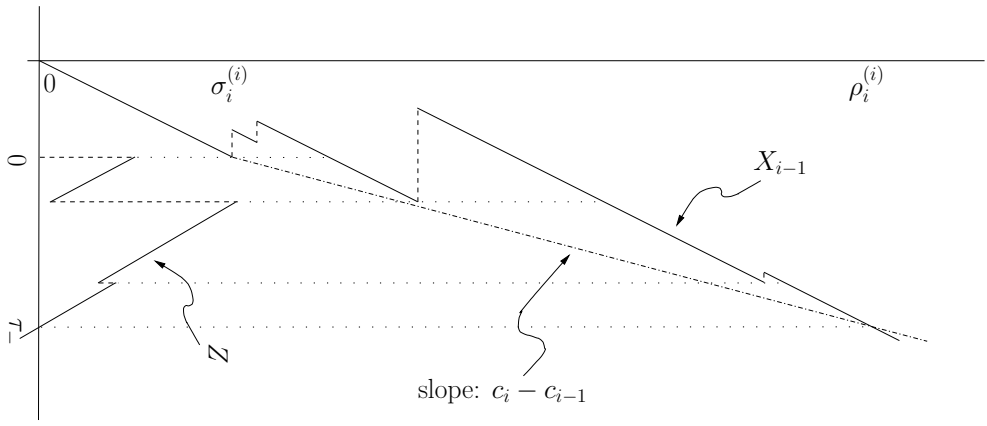}}
\caption{\label{fig:Zprocess}
Excursions of $\{X_{i-1}(t)- \underline X_{i-1}(t):t\ge0\}$
correspond to jumps of $Z$.}
\end{figure}

The idea is to apply Proposition~\ref{prop:excdistr} to the process
\[
Z(x):=
\inf\left\{t\geq 0: X_{i-1}\left(\sigma^{(i)}_i\right)-X_{i-1}\left(
\sigma^{(i)}_i+t\right)=x\right\}
-\frac x{c_{i-1}-c_i},
\]
see Figure~\ref{fig:Zprocess}.
In this diagram, excursions of $\{X_{i-1}(t)- \underline X_{i-1}(t):t\ge0\}$
correspond to jumps of $Z$. The relevant information
on the subexcursions is incorporated into $Z$ as jump marks.

Observe that $Z$ is a compound Poisson
process with negative drift $1/c_{i-1}-1/(c_{i-1}-c_i)$ and intensity $\lambda/c_{i-1}$,
starting with a (marked) jump at zero.
The jumps of $Z$ correspond to $(i-1)$-excursions, and the above marks are
assigned to the each of the jumps.
In terms of Proposition~\ref{prop:excdistr}, it remains to observe that
$\rho^{(i)}_{i}-\rho^{(i)}_{i-1}$ and
$\rho^{(i)}_i-\sigma^{(i)}_i$ correspond to $(\tau_--T_{N_-})/c_{i-1}$
and $\tau_-/(c_{i-1}-c_i)$ respectively.
\endproof

\medskip
With the recursion of Lemma~\ref{lem:recursionidle}, we can find the
joint transform of $\rho^{(k)}_k-\rho^{(k)}_j$ for $j=1,\ldots, k-1$,
which is required to work out Proposition~\ref{prop:analog}.
This is done in (\ref{eq:prdownkH}) below.
It is equivalent to find the transform of
$\rho^{(k)}_{j+1}-\rho^{(k)}_j$ for $j=1,\ldots,k-1$,
which is the content of the next proposition.
We have also added $\rho^{(k)}_k-\sigma^{(k)}_k$ for convenience.
The resulting formula has some remarkable features similar to
the formula in Theorem~\ref{thm:resultindep}.
Most interestingly, a quasi-product form appears here as well.

For $\beta\in\R_+^{k-1} \geq 0$, and $j=1,\ldots,k-1$, we define
\[
\C_j^k(\beta):=
c_j\sum_{\ell=j}^{k-1} \left(\frac1{c_{\ell+1}}-\frac1{c_\ell}\right)
(\lambda+\beta_\ell). 
\]

\begin{Prop}
\label{prop:distrho}
Suppose that $X$ is an $n$-dimensional L\'evy process
satisfying {\bf D} and {\bf H}.
Then for any $k= 2,\ldots,n$, $\beta\in\R_+^{k-1}$, $\gamma\geq 0$, we have
\begin{eqnarray*}
\mean \e{-\sum_{j=1}^{k-1} \beta_j \left(\rho_{j+1}^{(k)}-\rho^{(k)}_j\right)-
\gamma\left(\rho_k^{(k)}-\sigma_k^{(k)}\right)}&=&
\prod_{j=1}^{k-1} \frac{\beta_j+\lambda\mean \e{-\left[\C_{j+1}^k(\beta)+
\frac{c_{j+1}}{c_k}\gamma\right]
\left(\rho_{j+1}^{(j+1)}-\sigma^{(j+1)}_{j+1}\right)}}{\beta_j+\lambda
\mean \e{-\left[\C_{j}^k(\beta)+\frac{c_j}{c_k}\gamma\right]
\left(\rho_{j}^{(j)}-\sigma^{(j)}_{j}\right)}}\\
&&\mbox{}\times
\mean \e{-\left[\C_1^k(\beta)+\frac{c_1}{c_k}\gamma\right]\left(\rho_1^{(1)}-\sigma^{(1)}_1\right)}.
\end{eqnarray*}
\end{Prop}
\proof{Proof}
Since for $\ell=2,\ldots, i$, by definition of $\C_\ell^k(\beta)$,
\[
\left(\frac{c_{\ell-1}}{c_\ell}-1\right)(\lambda+\beta_{\ell-1}) + \frac{c_{\ell-1}}{c_\ell}
\C_\ell^k(\beta) = \C_{\ell-1}^k(\beta),
\]
it follows from Lemma~\ref{lem:recursionidle} that
\begin{eqnarray*}
\frac{\mean \e{-\sum_{j=1}^{\ell-1} \beta_j \left(\rho^{(\ell)}_{j+1}-\rho^{(\ell)}_j\right)
- \left[\C_\ell^k(\beta)+\frac{c_\ell}{c_k}\gamma\right]\left(\rho^{(\ell)}_\ell-\sigma^{(\ell)}_\ell\right)}}
{\mean \e{-\sum_{j=1}^{\ell-2} \beta_j \left(\rho^{(\ell-1)}_{j+1}-\rho^{(\ell-1)}_j\right)
-\left[\C_{\ell-1}^k(\beta)+\frac{c_{\ell-1}}{c_k}\gamma\right]
\left(\rho^{(\ell-1)}_{\ell-1}-\sigma^{(\ell-1)}_{\ell-1}\right)}}
&=&\frac{\beta_{\ell-1}+\lambda \mean \e{-\C_\ell^k(\beta,\gamma)
\left(\rho^{(\ell)}_\ell-\sigma^{(\ell)}_\ell\right)}}
{\beta_{\ell-1}+\lambda \mean \e{-\C_{\ell-1}^k(\beta,\gamma)
\left(\rho^{(\ell-1)}_{\ell-1}-\sigma^{(\ell-1)}_{\ell-1}\right)}}.
\end{eqnarray*}
The claim follows from this recursion (start with $\ell=k$ and note that $\C^k_k(\beta)=0$).
\endproof

\section{Multidimensional Skorokhod problems.}
\label{sec:skorokhod}
In the next sections, we apply results of the previous sections to
the analysis of {\em fluid networks}.
Such networks are closely related to (multidimensional) {\em
Skorokhod reflection problems}, which we describe first.
Subject to certain assumptions, we explicitly solve such
a reflection problem in Section~\ref{sec:specialSkorokhod}.
Section~\ref{sec:treenetworks} describes the fluid networks
associated to these special Skorokhod problems.

\medskip
Let $P$ be a nonnegative matrix with spectral radius
strictly smaller than one.
To a given c\`adl\`ag function $Y$ with values in $\R^n$ such that
$Y(0)=0$, one can associate a
c\`adl\`ag pair $(W,L)$ 
with the following properties ($w\in \R^n_+$):
\begin{description}
\item{\bf S1} $W(t) = w+Y(t) + (I-P')L(t), t\geq 0$,
\item{\bf S2} $W(t)\geq 0, t\geq 0$ and $W(0)=w$,
\item{\bf S3} $L(0)=0$ and $L$ is nondecreasing, and
\item{\bf S4} $\sum_{j=1}^n \int_0^\infty W_j(t)\, dL_j(t) = 0$.
\end{description}

It is known that such a pair exists and that it is unique; see Harrison
and Reiman~\cite{harrisonreiman:orthant1981} for the continuous case,
Robert~\cite{robert:stochnetworks2003} or Whitt~\cite[Thm.~14.2.3]{whitt:process2002}
for the c\`adl\`ag case,
and Kella~\cite{kella:reflecting2006} for a more general result.

It is said that $(W,L)$ is the solution to the Skorokhod problem of $Y$ in $\R^n_+$ with
reflection matrix $I-P'$ and initial condition $w$.

In general, the pair $(W,L)$ cannot be expressed explicitly in terms
of the driving process $Y$, with the notable exception of the
one-dimensional case. However, if the Skorokhod problem has a special structure,
this property carries over to a multidimensional setting.

\subsection{A special Skorokhod problem.}
\label{sec:specialSkorokhod}
It is the aim of this subsection to solve the Skorokhod problem for the pair $(W,L)$
under the following assumptions:
\begin{description}
\item{\bf N1} $P$ is strictly upper triangular,
\item{\bf N2} the $j$-th column of $P$ contains exactly one strictly positive element
for $j=2,\ldots,n$, and
\item{\bf N3} $Y_j$ is nondecreasing for $j=2,\ldots,n$.
\end{description}

In Section~\ref{sec:treenetworks}, we show that these assumptions
impose a `tree' structure on fluid networks.

\begin{Th}
\label{thm:skorokhodspecial}
Under {\bf N1--N3}, the solution to the Skorokhod problem of $Y$ in $\R_+^n$ is given by
\begin{eqnarray*}
L(t) &=& 0\vee\sup_{0\leq s\leq t} \left[-(I-P')^{-1} Y(s)-
(I-P')^{-1}w\right], \\
W(t) &=& w+Y(t)+ (I-P')L(t),
\end{eqnarray*}
where the supremum should be interpreted componentwise.
\end{Th}
\proof{Proof}
As $W$ is determined by $L$ and {\bf S1}, we only have to prove the expression
for $L$.
By Theorem~D.3 of Robert~\cite{robert:stochnetworks2003}, we know
that $L_i$ satisfies the fixed-point equation
\begin{equation}
\label{eq:toshowSkorokhod}
L_i(t)=0\vee\sup_{0\leq s\leq t} \left[(P' L)_i(s) - w_i-Y_i(s)\right]
\end{equation}
for $i=1,\ldots,n$ and $t\geq 0$.

As a consequence of {\bf N1}, we have $(I-P')^{-1}= I+P' +\ldots+P'^{n-1}$,
and the $j$-th row of $(I-P')^{-1}$ is the $j$-th row of
$I+P'+P'^2 +\ldots+P'^{j-1}$.
Therefore, the theorem asserts that
\begin{equation}
\label{eq:assertionSkorokhod}
L_i(t)=0\vee\sup_{0\leq s\leq t} \left[-\sum_{k=0}^{i-1}
\left[P'^{k}Y(s) +P'^{k}w)\right]_i\right].
\end{equation}
The proof goes by induction.
For $i=1$, (\ref{eq:assertionSkorokhod}) is the same equation as
(\ref{eq:toshowSkorokhod}).
Let us now suppose that we know that (\ref{eq:assertionSkorokhod})
holds for $i=1,\ldots,j-1$, where $j=2,\ldots, n$.
Furthermore, let $j^*<j$ be such that $p_{j^*j}>0$;
it is unique by {\bf N2}.
Equation (\ref{eq:toshowSkorokhod}) shows that
\begin{eqnarray}
L_{j}(t)&=&0\vee\sup_{0\leq s\leq t} \left[p_{j^*j} L_{j^*}(s) - w_j-Y_j(s)\right]\nonumber\\
&=&
0\vee\sup_{0\leq s\leq t} \left[\left(0\vee \sup_{0\leq u\leq s} -\sum_{k=0}^{j^*-1}
p_{j^*j} \left[P'^{k}Y(u)+P'^{k}w\right]_{j^*}\right) - w_j-Y_j(s)\right]\nonumber\\
&=& 0\vee\sup_{0\leq s\leq t} \left[\sup_{0\leq u\leq s} -\sum_{k=0}^{j^*-1}
p_{j^*j} \left[P'^{k}Y(u)+P'^{k}w\right]_{j^*} - w_j-Y_j(s)\right]
\label{from.n3.1}
\\ &=&
0\vee\sup_{0\leq u\leq t}\sup_{u\leq s\leq t}
\left[ -\sum_{k=0}^{j^*-1}\left[P'^{k+1} Y(u)+P'^{k+1}w\right]_j - w_j-Y_j(s)\right]\nonumber\\
&=&
0\vee\sup_{0\leq u\leq t}\left[ -\sum_{k=0}^{j^*}\left[P'^{k} Y(u)+P'^{k} w\right]_j\right],\label{from.n3.2}
\end{eqnarray}
where we have used {\bf N3} for the equalities \refs{from.n3.1} and \refs{from.n3.2}.

The proof is completed after noting that the $j$-th row
of $P'^k$ only contains zeroes for $k=j^*+1,\ldots, j-1$.~\endproof

\medskip
Instead of working directly with $W$, it is often convenient to
work with a transformed version, $\widetilde W:=(I-P')^{-1}W$.
The process $\widetilde{W}$ lies in a cone $\calC$, which is a polyhedron and a proper
subset of the orthant $\R^n_+$. Under the present
assumptions, at least one edge of $\calC$ is in the interior of $\R_+^n$ and at
least one is an axis. Below we give an interpretation of $\widetilde W$.

We next establish a correspondence between
the event that $W_j(t)=0$ and $\widetilde W_j(t)=0$
under an additional condition.

\begin{Prop}\label{p.equiv}
Suppose that {\bf N1--N3} hold, but with `nondecreasing' replaced by
`strictly increasing' in {\bf N3}.
Then we have $W_j(t)=0$ if and
only if $\widetilde W_j(t)=0$, for any $j=1,\ldots,n$
and $t\geq 0$.
\end{Prop}
\proof{Proof}
For $j=1$ we have $W_j(t)=\widetilde{W}_j(t)$, so the stated is satisfied;
suppose therefore that $j>1$.
Since the matrix $(I-P')^{-1}$ is lower triangular and nonnegative,
we straightforwardly get that
$\widetilde{W}_j(t)=0$ implies ${W}_j(t)=0$.

For the converse, observe that under {\bf N1--N2}
(see the proof of Theorem~\ref{thm:skorokhodspecial}; we use the same notation)
\[
\widetilde W_j(t) = \sum_{k=0}^{j-1} \left[P'^kW\right]_j(t)=
\sum_{k=0}^{j^*} \left[P'^kW\right]_j(t).
\]
An induction argument shows that it suffices to prove that
$W_j(t)=0$ implies $W_{j^*}(t)=0$.
To see that this holds, we observe that by {\bf S1} and
(\ref{eq:toshowSkorokhod}), $W_j(t)=0$ is equivalent to
\[
p_{j^*j }L_{j^*}(t)-w_j-Y_j(t)=0\vee \sup_{0\leq s\leq t}
\left[p_{j^*j}L_{j^*}(s)-w_j-Y_j(s)\right].
\]
The right-hand side of this equality is clearly nondecreasing.
Therefore, since $Y_j$ is strictly increasing by assumption, we conclude that
$dL_{j^*}(t)>0$, which immediately yields
$W_{j^*}(t)=0$ by {\bf S4}. This completes the proof.
\endproof

\subsection{L\'evy-driven tree fluid networks.}
\label{sec:treenetworks}
In this subsection, we define a class of L\'evy-driven
fluid networks, which we call {\em tree fluid networks}.
We are interested in the steady-state behavior
of such networks.

Consider $n$ (infinite-buffer) fluid queues,
with external input to queue $j$ in the time interval
$[0,t]$ given by $J_j(t)$.
We assume that
$J=\{J(t):t\geq 0\}=\{(J_1(t),\ldots,J_n(t))^{'}:t\geq 0\}$ is a
c\`adl\`ag L\'evy process
starting in $J(0)= 0\in \R^n_+$.
The buffers are continuously drained at a constant rate
as long as there is content in the buffer.
These drain rates are given by a vector $r$; for buffer $j$, the rate is $r_j>0$.

The interaction between the queues is modeled as follows.
A fraction $p_{ij}$ of the output of station $i$ is immediately
transferred to station $j$, while a fraction $1-\sum_{j\neq i} p_{ij}$ leaves the system. We
set $p_{ii}=0$ for all $i$, and suppose that $\sum_j p_{ij}\le 1$.
The matrix $P=\{p_{ij}:i,j=1,\ldots,n\}$ is called the {\em
routing matrix}.
We assume that for any station $i$, there is at most one
station feeding buffer $i$, and that $p_{ij}=0$ for $j<i$.
The resulting network can be represented by a (directed) tree.
Indeed, the stations then correspond to nodes, and there is a vertex
from station $i$ and $j$ if $p_{ij}>0$.
This motivates the name `tree fluid networks'.
We represent such a fluid network by the triplet $(J,r,P)$.
Note that $P$ satisfies {\bf N1--N2} by definition
of a tree fluid network.

The {\em buffer content} process $W$ and {\em regulator} $L$ associated
to the fluid network $(J,r,P)$ are defined as the solution of the
Skorokhod problem of
\[
Y(t):=J(t) - (I-P')rt
\]
with reflection matrix $I-P'$.
The buffer content is sometimes called the {\em workload},
explaining the notation $W$.
Importantly, the dynamics of the network are given by {\bf S1--S4}, as the
reader may verify.
The process $L_j$ can be interpreted as the cumulative unused capacity
in station $j$.

Associated to the processes $W$ and $L$, one can also define the
process of the {\em age of the busy period}: for $j=1,\ldots,n$, we set
\begin{equation}
\label{eq:defbusyperiod}
B_j(t):=t-\sup\{s\leq t: W_j(s)=0\},
\end{equation}
and let $B(t)=(B_1(t),\ldots,B_n(t))'$.
Hence, if there is work in queue $j$ at time $t$ (that is, $W_j(t)>0$),
$B_j(t)$ is the
time that elapsed after the last time that the $j$-th queue
was empty.
If there is no work in queue $i$ at time $t$, then $B_i(t)=0$.
Similarly, one can also define the {\em age of the idle period}
for $j=1,\ldots,n$:
\[
I_j(t):=t-\sup\{s\leq t: W_j(s)\neq 0\},
\]
and the corresponding vector $I(t)$.
As a result of these definitions, $I_j(t)>0$ implies $B_j(t)=0$ and
$B_j(t)>0$ implies $I_j(t)=0$ for $j=1,\ldots,n$.
The quantities $\widetilde B_j(t)$ and $\widetilde I_j(t)$
are defined similarly,
but with $W_j$ replaced by the $j$-th element of $\widetilde W=(I-P')^{-1}W$.

The random variables $\widetilde W_j$, $\widetilde B_j$, and $\widetilde I_j$
have a natural interpretation.
Indeed, let us consider all stations on a path from the root of the tree
to station $j$. The {\em total} content of the buffers along this path is then given by
$\widetilde W_j$. Consequently, $\widetilde B_j$ and $\widetilde I_j$ correspond to the
ages of the busy and idle periods of this aggregate buffer.


In the rest of the paper, we
assume that the tree fluid network has the following additional
properties:
\begin{description}
\item[{\bf T1}] If $p_{ij}>0$, then $p_{ij}> r_j/r_i$,
\item[{\bf T2}] $J_j(t)$ are nondecreasing for $j=2,\ldots,n$,
\item[{\bf T3}] $J$ is an $n$-dimensional L\'evy process, and
\item[{\bf T4}] $J$ is integrable and $(I-P')^{-1} \mean J(1)<r$.
\end{description}
An important consequence of {\bf T1} and {\bf T2} is that $Y$ is
componentwise nondecreasing, except for $Y_1$. Consequently, if
{\bf T1} and {\bf T2} hold for a tree fluid network, then {\bf
N1--N3} are automatically satisfied for the associated Skorokhod
problem. Hence, Theorem~\ref{thm:skorokhodspecial} gives an
explicit description of the buffer contents in the network.
Note that {\bf T4} ensures {\em stability} of the network.

Let us now define the process
\[
X(t):=(I-P')^{-1} Y(t) =(I-P')^{-1} J(t)-rt.
\]
In view of assumption {\bf T1},
the down-stream buffer contents always grow when one of the up-stream buffers
is nonempty.
Moreover, under {\bf T1},
$\widetilde W$ is itself a reflected process, that is
$(\widetilde {W},\widetilde {L})$ is the solution to the Skorokhod problem for
$X$ with reflection matrix $I$ and initial condition $(I-P')^{-1}w$.
Therefore, each coordinate of $\tilde W$ is a one-dimensional reflected process.
A similar assumption facilitates the analysis in \cite{kella:paralleldependent1993},
\cite[Thm.~4.1 and Lem.~4.2]{kella:stationarityfeedfw1997},
and \cite{kellawhitt:tandem1992}.

In the next proposition, we find the steady-state behavior
of the buffer content and the age of the busy (and idle) period for
the L\'evy-driven tree fluid network $(J,r,P)$.
We also consider the case where the inequality $p_{ij}>r_j/r_i$ in {\bf T1} holds
only weakly (i.e. $p_{ij}\ge r_j/r_i$),
as this plays a role in priority fluid systems (see Section~\ref{sec:priority} below).

Recall the definitions of $G=G^X$ and $H=H^X$
in Sections~\ref{sec:splitmax} and \ref{sec:splitlastexit} respectively.

\begin{Prop}
\label{prop:skorokhod}
Suppose that {\bf T1--T4} hold for the tree fluid network $(J,r,P)$.
\begin{itemize}
\item[(i)] For any initial condition $W(0)=w$,
the triplet of vectors $(W(t), B(t), I(t))$ converges in
distribution to $((I-P') \overline{X}, G^X, H^X)$ as $t\to\infty$.
\item[(ii)] If the second inequality in {\bf T1} holds only weakly,
then for any initial condition $W(0)=w$,
the triplet of vectors $(W(t), \widetilde B(t), \widetilde I(t))$ converges in
distribution to $((I-P') \overline{X}, G^X, H^X)$ as $t\to\infty$.
\end{itemize}
\end{Prop}
\proof{Proof}
Throughout this proof, a system of equations like (\ref{eq:defbusyperiod})
is abbreviated by $B(t)=t-\sup\{s\leq t: W(s)=0\}$.

We start with the proof of (ii).
By Theorem~\ref{thm:skorokhodspecial}, we have for
any $t>0$
\[
\widetilde W(t) = [x+X(t)]\vee \sup_{0\leq s\leq t}[X(t)-X(s)],
\]
where $x=(I-P')^{-1}w$.
Moreover, as a consequence of Proposition~\ref{p.equiv}, we have
\begin{eqnarray*}
\widetilde B(t)&=&t-\sup\{s\leq t:\widetilde W(s)=0\}\\
&=&t-\sup\left\{s\leq t:
x+X(s)=0\wedge \inf_{0\leq u\leq s} [x+X(u)]\right\}\\
&=& t-\sup\left\{s\leq t:
x+X(s)=0\wedge \inf_{0\leq u\leq t} [x+X(u)]\right\},
\end{eqnarray*}
where the last equality is best understood by sketching a sample path of $X$.
The supremum over an empty set should be interpreted as zero.

This reasoning carries over to idle periods:
\[
\widetilde I(t)=t-\sup\left\{s\leq t:
x+X(s)\neq 0\wedge \inf_{0\leq u\leq s} [x+X(u)]\right\}.
\]

Due to the stationarity of the increments of $\{X(t), \ t\ge 0\}$
({\bf T3}), we may
extend $X$ to the two-sided process $\{X(t), \ t\in\R\}$.
This leads to
\[
\left(\begin{array}{c}\widetilde{W}(t)\\ \widetilde B(t)\\\widetilde I(t)
\end{array}\right) =_{\rm d} \left(\begin{array}{c}
[x-X(-t)]\vee\sup_{-t\le s\le 0}[-X(s)]\\
-\sup\left\{s: -t\le s\le 0,
-X(s)= [x-X(-t)]\vee \sup_{-t\le u\le 0}[-X(u)]\right\}\\
-\sup\left\{s: -t\le s\le 0,
-X(s)\neq [x-X(-t)]\vee \sup_{-t\le u\le s}[-X(u)]\right\}
\end{array}\right).
\]
Since $x-X(-t)\to-\infty$ almost surely by {\bf T4},
this tends to
\[
\left(\begin{array}{c}
\sup_{s\le 0}[-X(s)]\\
-\sup\left\{s\le 0: -X(s)=\sup_{u\le 0}[-X(u)]\right\}\\
-\sup\left\{s\le 0: -X(s)\neq \sup_{u\le s}[-X(u)]\right\}
\end{array}\right),
\]
a vector that is almost surely finite, again by {\bf T4}.
By time-reversibility
(see Lemma~II.2 of Bertoin~\cite{bertoin:levy1996}),
the latter vector is equal in distribution to $(\overline X, G^X, H^X)$.

The first claim follows from (ii) after noting that
$B(t)=\widetilde B(t)$ and $I(t)=\widetilde I(t)$ by
Proposition \ref{p.equiv}.
\endproof

\medskip
We remark that the above proof does not use {\bf T3} to the fullest.
Indeed, for the proposition to hold, it suffices that $J$ has stationary
increments and that it is time-reversible.

\medskip
Let us now suppose that the initial buffer content $w$ is random.
Proposition~\ref{prop:skorokhod} shows, after a standard argument,
that $\{W(t)\}$ is a stationary process if $W(0)=w$ is distributed
as $\mu^*$, where $\mu^*$ is the distribution of $(I-P')\overline X$.
We now show that this stationary distribution is unique.

\begin{Cor}
\label{cor:uniquenessskorokhod}
Suppose that {\bf T1--T4} hold for the tree fluid network $(J,r,P)$.
Then $\mu^*$ is the only stationary distribution.
\end{Cor}
\proof{Proof}
Suppose there exists another stationary distribution $\mu^*_0\neq \mu^*$.
Let $W^*_0$ be the corresponding stationary process.
For any Borel set $B$ in $\R_+^n$ and any $t\geq 0$, we then have
$
\pr (W^*_0(0)\in B)=\pr (W^*_0(t)\in B)
$.
Therefore,
\begin{eqnarray}
\pr (W^*_0(0)\in B)
&=&
\lim_{t\to\infty}\pr (W^*_0(t)\in B)\nonumber\\
&=&
\lim_{t\to\infty} \int_0^\infty \pr(W^*_0(t)\in B|W^*_0(0)=w)\pr (W^*_0(0)\in dw)\nonumber\\
&=&
 \int_0^\infty \lim_{t\to\infty}\pr(W^*_0(t)\in B|W^*_0(0)=w)\pr (W^*_0(0)\in dw)\nonumber\\
&=&
 \int_0^\infty \pr((I-P')  \overline{X}\in B) \pr (W^*_0(0) \in dw)= \pr((I-P')
\overline{X}\in B),\nonumber
\end{eqnarray}
where the second last equation is due to Proposition~\ref{prop:skorokhod}.
This is clearly a contradiction.
\endproof

\medskip
Corollary~\ref{cor:uniquenessskorokhod} answers, for the special case of tree fluid networks,
a question from the paper of
Konstantopolous, Last and Lin~\cite{konstantopoulos:networks2004}
on the uniqueness of the stationary distribution.
Note that for the queueing problem related to $(J,r,P)$,
the uniqueness of the stationary
distribution was discussed in
Kella~\cite{kella:stationarityfeedfw1997}.
In contrast to the setting in \cite{kella:stationarityfeedfw1997},
we allow for the first component of $J(t)$ to be a general L\'evy process.

In the next section,
we combine Proposition~\ref{prop:skorokhod} with
the results given in Sections~\ref{sec:distrXG} and \ref{sec:distrH} to
study particular networks.


\section{Tandem networks and priority systems.}
\label{sec:tandem}
In this section, we analyze $n$ fluid queues in tandem,
which is a tree fluid network with a special structure.
We also analyze a closely related priority system.

The tandem structure is specified by the form of the routing matrix:
we suppose that $P$ is such that
$p_{i,i+1}>0$ for $i=1,\ldots,n-1$, and $p_{ij}=0$ otherwise.
Observe that we allow $p_{i,i+1}>1$,
and that it is not really a restriction to exclude $p_{i,i+1}=0$;
otherwise the queueing system splits into independent tandem
networks.

In all of our results, we suppose that the tandem
system $(J,r,P)$ satisfies {\bf T1--T4}.
We rule out the degenerate case where the first $j\geq 1$
components of $J$ are deterministic drifts,
since an equivalent problem can then be studied with the first $j$ stations
removed.
We also impose the following assumptions on the input L\'evy
process $J$:

\begin{description}
\item{\bf T5} $J$ has mutually independent components, and
\item{\bf T6} The L\'evy measure of $J_1$ is supported on $\R_+$.
\end{description}

Observe that under {\bf T2--T3}, {\bf T5} implies that
$J_2,\ldots,J_n$ are independent nonnegative subordinators.

\medskip
This section consists of three parts.
In Section~\ref{sec:tandemgeneralities},
we are interested in the joint (steady-state) distribution
of the buffer contents and the ages of the busy periods
for fluid tandem networks, i.e.,
in the distribution of $(W(\infty), B(\infty))$.
Section~\ref{sec:singlecompound} considers the situation
of a single compound Poisson input to the system.
For that system, we are also interested in the
ages of the idle periods, i.e., in the
vector $I(\infty)$.
In Section \ref{sec:priority}, we analyze buffer contents and busy
periods in a priority system.

\subsection{Generalities.}
\label{sec:tandemgeneralities} To find the joint distribution of
$W(\infty)$ and $B(\infty)$, throughout this section denoted by $W$
and $B$ respectively, we rely on Proposition~\ref{prop:skorokhod}.
This motivates the analysis of $X(t)=(I-P')^{-1}J(t)-rt$. For
$i=2,\ldots,n$, we define the cumulant of $J_i(t)$ by
$\theta_i^J(\beta):=-\log\mean\e{-\beta J_i(1)}$, $\beta\geq
0$. As in Section~\ref{sec:distrH}, we
write $\psi_i$ (defined by $\psi_i(\beta)=\log\mean\e{-\beta
X_i(1)}$) for the Laplace exponent of $-X_i$. Its inverse is again
denoted by $\Phi_i$.

Under {\bf T2} and {\bf T6},
the L\'evy measure of $X$ is
supported on $\R_+^n$. Moreover, as we ruled out trivial queues in the network,
each of the components of $\overline X$ has a nondegenerate distribution.
Therefore, let us recall that the following holds
(see, e.g., Theorem VII.4 in \cite{bertoin:levy1996}):
for $\alpha,\beta\geq 0$, $(\alpha,\beta)\neq (0,0)$, $\beta\neq \Phi_i(\alpha)$,
$i=1,\ldots,n$, we have
\begin{equation}
\mean \e{-\alpha G_i-\beta \overline X_i }= -\mean
X_i(1)\frac{\Phi_i\left(\alpha \right)-\beta}
{\alpha -\psi_i(\beta) }.
\label{eq:fluctuation}
\end{equation}
This identity plays a crucial role in the results of this
section. For notational convenience, we shall write that
(\ref{eq:fluctuation}) holds for any $\alpha,\beta \geq 0$,
without the requirements $(\alpha,\beta)\neq (0,0)$ and $\beta \neq \Phi_i(\alpha)$.

Now we can formulate the main result of this subsection.
We remark that the first formula also holds if $J_1$ is
not necessarily spectrally positive.
For instance, it allows for phase-type downward jumps;
see \cite{dieker:fluclevy2006} for the joint transform
of $\oX_j$ and $G_j$ in that case.

\begin{Th}
\label{th:tandem}
Consider a tandem fluid network $(J,r,P)$ for which {\bf T1--T6}
holds.
Then for $\omega,\beta\in \R_+^n$, the transform
$\mean\e{-\langle\omega,W\rangle- \langle\beta,B\rangle}$
equals
\[
\prod_{j=1}^{n-1}
\frac{\mean \e{ -\left[\sum_{\ell=j+1}^n
\theta^J_{\ell}\left(\omega_\ell\right)+
 \sum_{\ell=j+1}^n   (p_{\ell-1\ell}r_{\ell-1}-r_{\ell})\omega_\ell+
\sum_{\ell=j}^n \beta_\ell \right] G_{j}-\omega_j \overline X_{j}}} {\mean
\e{ -\left[\sum_{\ell=j+1}^n
\theta^J_{\ell}\left(\omega_\ell\right)+
 \sum_{\ell=j+1}^n   (p_{\ell-1\ell}r_{\ell-1}-r_{\ell})\omega_\ell+
\sum_{\ell=j+1}^n \beta_\ell \right] G_{j}-p_{j,j+1}\omega_{j+1} \overline X_{j}}}
\times \mean \e{-\beta_n G_{n}-\omega_n
\overline X_{n}}.
\]
Consequently, we have for $\omega,\beta\in\R_+^n$,
\begin{eqnarray*}
\lefteqn{ \mean\e{-\langle\omega,W\rangle- \langle\beta,B\rangle}=
       -\mean X_n(1)\frac{\Phi_n\left(\beta_n \right)-\omega_n}
       { \beta_n -\psi_n(\omega_n) }
}\\
&\times& \prod_{j=1}^{n-1}
  \frac{\Phi_j\left(\sum_{\ell=j+1}^n \theta^J_{\ell}\left(\omega_\ell\right)+
         \sum_{\ell=j+1}^n   (p_{\ell-1,\ell}r_{\ell-1}-r_{\ell})\omega_\ell+
         \sum_{\ell=j}^n \beta_\ell \right)-\omega_j}
       { \Phi_j\left(\sum_{\ell=j+1}^n \theta^J_{\ell}\left(\omega_\ell\right)+
         \sum_{\ell=j+1}^n   (p_{\ell-1,\ell}r_{\ell-1}-r_{\ell})\omega_\ell+
         \sum_{\ell=j+1}^n \beta_\ell \right)-p_{j,j+1}\omega_{j+1} }\\
 &\times& \prod_{j=1}^{n-1}
  \frac{\sum_{\ell=j+1}^n \theta^J_{\ell}\left(\omega_\ell\right)+
         \sum_{\ell=j+1}^n   (p_{\ell-1,\ell}r_{\ell-1}-r_{\ell})\omega_\ell+
         \sum_{\ell=j+1}^n \beta_\ell-\psi_j(p_{j,j+1}\omega_{j+1}) }
       {\sum_{\ell=j+1}^n \theta^J_{\ell}\left(\omega_\ell\right)+
         \sum_{\ell=j+1}^n   (p_{\ell-1,\ell}r_{\ell-1}-r_{\ell})\omega_\ell+
         \sum_{\ell=j}^n \beta_\ell-\psi_j(\omega_j)}.
\end{eqnarray*}

\end{Th}
\proof{Proof}
By Proposition~\ref{prop:skorokhod}(i),
$(W,B)=_{\rm d}((I-P')\overline X, G^X)$.
Hence we have
\begin{equation}
\mean\e{-\langle\omega,W\rangle-\langle\beta,B\rangle} =
\mean\e{-\langle(I-P)\omega,(I-P')^{-1}W\rangle-\langle\beta,B\rangle}
= \mean\e{-\langle\beta,G\rangle-\langle(I-P)\omega,\overline
X\rangle}.\label{from.skorokhod}
\end{equation}

Now note that the stability condition {\bf T4} for $(J,r,P)$ implies
{\bf D} for $X$ by the law of large numbers.
Thus, in order to apply Theorem
\ref{thm:resultindep} for \refs{from.skorokhod}, it is enough to
check that {\bf G} holds. Standard algebraic
manipulations give
\[
X_1(t)=J_1(t)-r_1t
\]
and
\[
X_{i+1}(t)=p_{i,i+1} X_i(t)+J_{i+1}(t)+(p_{i,i+1}r_i-r_{i+1})t
\]
for $i=1,\ldots,n-1$. Hence, {\bf G} holds with
$K_i=p_{i-1,i}$ and
$\Upsilon_i(t)=J_{i}(t)+(p_{i-1,i}r_{i-1}-r_{i})t$.

As a result, we know that from
Theorem \ref{thm:resultindep},
\begin{eqnarray*}
\lefteqn{ \mean\e{-\langle\beta,G\rangle-\langle(I-P)\omega,\overline X\rangle} =
\mean\e{-\langle\beta,G\rangle-\langle\widetilde{\omega},\overline X\rangle}}
\\
&=& \prod_{j=1}^{n-1} \frac{\mean \e{
-\left[\sum_{\ell=j+1}^n \theta^\Upsilon_{\ell} \left(\sum_{k=\ell}^n
K_\ell^k \widetilde{\omega}_k\right)+ \sum_{\ell=j}^n \beta_\ell \right]
G_{j}-\left(\sum_{k=j}^n K_j^k
\widetilde{\omega}_k\right) \overline X_{j}}}
{\mean \e{
-\left[\sum_{\ell=j+1}^n \theta^\Upsilon_{\ell} \left(\sum_{k=\ell}^n
K_\ell^k \widetilde{\omega}_k\right)+ \sum_{\ell=j+1}^n \beta_\ell
\right] G_{j}-\left(\sum_{k=j+1}^n K_j^k
\widetilde{\omega}_k\right) \overline X_{j}}}
\times \mean \e{-\beta_n G_{n}-\widetilde{\omega}_n \overline
X_{n}},
\end{eqnarray*}
where we have set $\widetilde{\omega}=(I-P)\omega$ for notational convenience.

The reader may check that
$\sum_{k=j}^n K_j^k \widetilde{\omega}_k
=\omega_j$
and
$\sum_{k=j+1}^n K_j^k \widetilde{\omega}_k=
p_{j,j+1}\omega_{j+1}$,
leading to the first claim.
The second assertion is a consequence of the first and
\refs{eq:fluctuation}.
\endproof

\medskip
Theorem~\ref{th:tandem} extends several results
from the literature on
the steady-state distribution of the buffer content for tandem L\'evy networks. In
particular, if
$J(t)=(J_1(t),0)'$, $P=(p_{ij})$, with $p_{12}=1$ and zeroes elsewhere,
if one chooses
$\beta_1=\beta_2=0$ and $\omega_1=0$ in Theorem \ref{th:tandem},
then one obtains Theorem~3.2 of D\c{e}bicki, Mandjes and van
Uitert~\cite{debicki:tandem2005}. Additionally, if one chooses
$\beta_1=\beta_2=0$ and supposes that $J_1$ is a subordinator,
we recover the results of
Kella~\cite{kella:paralleldependent1993}.

Even if the Laplace transform of $(G_j,\overline X_j)$ can be inverted, it
is generally not straightforward to invert the Laplace transform of $(W,B)$
given in Theorem~\ref{th:tandem}.
Some progress has been recently made in case $n=2$; for a Brownian fluid system,
Lieshout and Mandjes~\cite{lieshoutmandjes:brownian2006} calculate the
distribution of $W$.
Avram, Palmowski, and Pistorius~\cite{avram:twodimensional2006}
study a compound Poisson setting with exponential jumps.
A different type of explicit solution can be found in the work of
Harrison~\cite{harrison:tandem1978}; he gives an example closely
related to the framework of the present paper.

For use in Section~\ref{sec:priority}, we point out that the expression in
Theorem~\ref{th:tandem} is $\mean \e{-\langle\omega, W\rangle-
\langle \beta, \widetilde B\rangle}$ if the second inequality in
{\bf T1} is weak, cf.~Proposition~\ref{prop:skorokhod}(ii).

\subsubsection*{The lengths of the busy periods.}
Besides the Laplace transforms of the {\em ages} $B$ of the
busy periods, Theorem~\ref{th:tandem} also enables us to find the
Laplace transforms of the {\em length} $V$ of the
steady-state running busy periods. Indeed, let $D_i$,
$i=1,\ldots,n$ denote the steady-state remaining lengths of the
running busy period, so that $V_i=B_i+D_i$. We know that $D_i$ and $B_i$
are equal in distribution. In fact, following for instance \cite[Sec.~V.3]{asmussen:applprob2003},
we have
\begin{eqnarray}
(B_i,D_i)=_{\rm d}(U_i V_i, (1-U_i) V_i),\label{eq.bd}
\end{eqnarray}
where $U_i$ are i.i.d. and uniform on $[0,1]$.

For the Brownian (single-station) fluid queue, the
following result is Corollary~3.8 of
Salminen and Norros~\cite{salminennorros:brownianstorage2001}.

\begin{Cor}
Consider a tandem fluid network $(J,r,P)$ for which {\bf T1--T6}
holds. Then for $\alpha,\beta\geq 0$, $\alpha\neq \beta$,
\[
\mean \e{-\alpha B_i -\beta D_i} = -\mean
X_i(1)\frac{\Phi_i(\alpha)-\Phi_i(\beta)}{\alpha-\beta}.
\]
Moreover, we have for $\alpha\geq 0$,
\[
\mean \e{-\alpha V_i} =-\mean X_i(1)\frac{d
\Phi_i(\alpha)}{d \alpha}.
\]
\end{Cor}
\proof{Proof} Since the second claim follows straightforwardly from
the first, we only prove the first expression. Following (\ref{eq.bd}), we have
for $\alpha\neq \beta$,
\begin{eqnarray*}
(\alpha-\beta)\mean \e{-\alpha B_i -\beta D_i} &=&
(\alpha-\beta) \mean \e{-(\alpha-\beta)U_iV_i-\beta V_i}=
(\alpha-\beta)\mean \int_0^1 \e{-(\alpha-\beta)u V_i-\beta V_i}du\\
&=& \mean \int_\beta^\alpha \e{-u V_i} du
= \mean \int_0^\alpha \e{-u V_i} du-\mean \int_0^\beta \e{-u V_i} du.
\end{eqnarray*}
The two identities that result upon setting $\beta=0$
and $\alpha=0$ can be used to express the first and second
expectation in terms of the
Laplace transform of $B_i$ and $D_i$ respectively; this
yields for $\alpha\neq \beta$
\[
\mean \e{-\alpha B_i -\beta D_i} =\frac{1}{\alpha-\beta}
\left[\alpha\mean \e{-\alpha B_i} - \beta \mean \e{-\beta
B_i}\right],
\]
where we have used the equality in distribution of $B_i$ and $D_i$.
Application of \refs{eq:fluctuation} completes the proof.~\endproof

\subsection{A single compound Poisson input.}
\label{sec:singlecompound}
In this subsection,
we examine a tandem fluid network with a single compound Poisson input
\cite{kellawhitt:tandem1992}.
The following assumption formalizes our framework.

\begin{description}
\item[{\bf T7}] $p_{i,i+1}=1$ for $i=1,\ldots,n-1$, while $p_{ij}=0$ otherwise, and
\item[{\bf T8}] $J_1$ is a compound Poisson process with positive drift
$\mathtt d$ and intensity $\lambda$,
and $J_j\equiv0$ for $j=2,\ldots,n$.
Moreover, $r_j$ decreases strictly in $j$ and $\mean J(1)<r_n$.
\end{description}

An important consequence of {\bf T7} and {\bf T8} is that
\begin{equation}
\label{eq:simplification}
(r_j-r_k)\omega=\psi_j(\omega)  -\psi_k(\omega),
\end{equation}
which simplifies the resulting expressions
in view of fact that we often deal with ratios of the
fluctuation identity (\ref{eq:fluctuation}).
Interestingly, it is also possible
to study (joint distributions of) idle periods under these assumptions.

The following corollary collects some results that follow from
{\bf T7} and {\bf T8} and Theorem~\ref{th:tandem}.
Many
interesting formulas can be derived, but we have selected two examples for which the formulas are
especially appealing.

\begin{Cor}
\label{cor:tandem}
Consider a tandem fluid network $(J,r,P)$ for which {\bf T7--T8} holds.
\begin{enumerate}
\item[(i)] For $i=1,\ldots,n$, and $\omega,\beta\geq 0$,
we have
\[
\mean \e{-\omega W_i-\beta B_i}=
-\mean X_i(1) \frac{\Phi_i(\beta)-\omega}{\beta+(r_{i-1}-r_i)\omega}
\times \frac{\Phi_{i-1}((r_{i-1}-r_i)\omega+\beta)}
{\Phi_{i-1}((r_{i-1}-r_i)\omega+\beta)-\omega}.
\]
Moreover, $\pr(W_i=0)=\pr(B_i=0) = \frac{\mean X_i(1)}{\mathtt d-r_i}$.
\item[(ii)] For $i=2,\ldots,n$, $\omega,\beta\geq 0$, we have
\[
\mean \left[\e{-\omega W_{i} -\beta B_{i}};
W_{i-1}=0\right]=  -\frac{\mean X_{i}(1)}{\mathtt d-r_{i-1}}
\frac{\Phi_{i}(\beta)-\omega}
{\Phi_{i-1}((r_{i-1}-r_{i})\omega+\beta) -\omega}.
\]
\end{enumerate}
\end{Cor}
\proof{Proof}
To prove (i), apply Theorem~\ref{th:tandem} to obtain for $i=1,\ldots,n$,
\[
\mean \e{-\omega W_i -\beta B_i} =
\frac{\mean \e{-\left[(r_{i-1}-r_i)\omega
+\beta \right]G_{i-1}}}{\mean \e{-\left[(r_{i-1}-r_i)\omega
+\beta \right]G_{i-1}-\omega \oX_{i-1}}} \mean \e{-\beta G_i-\omega \oX_i}.
\]
With \refs{eq:fluctuation},
this leads immediately to the given formula
after invoking (\ref{eq:simplification}).

We find $\pr(W_i=0)$ upon choosing $\omega=0$ and noting that
\[
\pr(W_i=0)=\pr(B_i=0) = \lim_{\beta\to\infty} \mean \e{-\beta G_i} =
-\mean X_i(1) \lim_{\beta\to\infty} \frac{\Phi_i(\beta)}\beta =
\frac{\mean X_i(1)}{\mathtt d-r_i},
\]
where the last equality follows from Proposition I.2 in
\cite{bertoin:levy1996}.

The second claim uses a similar argument;
it follows from Theorem~\ref{th:tandem} that for $i=2,\ldots,n$
\[
\mean \e{-\omega_i W_i-\beta_{i-1} B_{i-1}-\beta_i B_i} =
\frac{\mean \e{- [(r_{i-1}-r_i)\omega_i +\beta_{i-1} +\beta_i] G_{i-1}}}
{\mean \e{-[(r_{i-1}-r_i)\omega_i +\beta_i] G_{i-1}-\omega_i \oX_{i-1}}}
\mean \e{-\beta_i G_i-\omega_i \oX_i},
\]
and the numerator of the fraction tends
to $\pr(W_{i-1}=0)$ as $\beta_{i-1}\to\infty$.
Now apply \refs{eq:fluctuation} and \refs{eq:simplification}.
\endproof

\medskip
We end this subsection with an application of the theory
in Section~\ref{sec:distrH}, which
enables us to study the idle periods in a tandem fluid network
satisfying {\bf T7--T8}.
For $\gamma\in\R^{k-1}_+$, we set
\[
\Dcal_j^k(\gamma):=c_j\sum_{\ell=j}^{k-1}\left(\frac{1}{c_{\ell+1}}-
\frac1{c_\ell}\right)\left(\lambda+\sum_{p=1}^\ell\gamma_p\right),
\]
which is similar to the definition of $\mathcal C_j^k$
in Section~\ref{sec:distrH}.

\begin{Prop}
Consider a tandem fluid network $(J,r,P)$ for which {\bf T7--T8} holds.
For $\gamma\in\R^n_+$, we have
\[
\mean \e{-\langle \gamma, I\rangle} = 1-\sum_{k=1}^n
\pr(W_k=0) \meandown_k\left[\e{-\sum_{\ell=1}^{k-1} \gamma_\ell H_\ell}
\left(1-\e{-\gamma_k H_k}\right)\right],
\]
where $\pr(W_j=0)$ is given in Corollary~\ref{cor:tandem}(i), and
\begin{eqnarray}
\meandown_k \e{-\sum_{\ell=1}^k \gamma_\ell H_\ell}&=&
\frac{\lambda+\sum_{\ell=1}^{k-1} \gamma_\ell\left(1-\frac{c_k}{c_\ell}\right) -c_k
\Phi_1(\Dcal_1^k(\gamma))}{\lambda+\sum_{\ell=1}^k \gamma_\ell}
\nonumber \\ &&\mbox{}\times
\prod_{j=1}^{k-1}
\frac{\lambda +\sum_{\ell=1}^{k-1}\gamma_\ell-\sum_{\ell=j+1}^{k-1} \frac{c_k}{c_{\ell}} \gamma_{\ell}
-c_k\Phi_{j+1}(\Dcal^k_{j+1}(\gamma))}{\lambda +\sum_{\ell=1}^{k-1}\gamma_\ell
-\sum_{\ell=j+1}^k \frac{c_k}{c_{\ell}}
\gamma_{\ell} -c_k\Phi_{j}(\Dcal^k_{j}(\gamma))}.
\label{eq:prdownkH}
\end{eqnarray}
\end{Prop}
\proof{Proof}
Note that {\bf T7} and {\bf T8} imply {\bf H}.
The first claim follows from Proposition~\ref{prop:skorokhod} and
the facts that for $k=2,\ldots,n$,
\[
\mean \e{-\sum_{\ell=1}^k\gamma_\ell H_\ell} = \mean \e{-\sum_{\ell=1}^{k-1}
\gamma_\ell H_\ell} +\meandown_k\left[\e{-\sum_{\ell=1}^{k-1}
\gamma_\ell H_\ell}\left(1-\gamma_k H_k\right)\right] \pr(\oX_k=0),
\]
and $\mean \e{\gamma_1 H_1} = 1- \meandown_1 \left[1-\e{-\gamma_1H_1}\right]
\pr(\oX_1=0)$. These identities follow after observing that $H_k$ vanishes on the event
$\{\oX_k=0\}$, and that $\{\oX_k=0\}$ is
the complement of $\{\oX_k>0\}$.

Let us now prove the expression for the $\prdown_k$-distribution
of $(H_1,\ldots,H_k)'$. From Proposition~\ref{prop:analog} and
Proposition~\ref{prop:distrho}, we know that
\[
\meandown_k \e{-\sum_{\ell=1}^k \gamma_\ell H_\ell}=
\frac{\lambda\mean \e{-\Dcal_1^k(\gamma)\left(\rho_1^{(1)}-\sigma^{(1)}_1\right)}}
{\lambda+\sum_{\ell=1}^k\gamma_\ell}
\prod_{j=1}^{k-1} \frac{\sum_{\ell=1}^j\gamma_\ell+\lambda\mean \e{-\Dcal_{j+1}^k(\gamma)
\left(\rho_{j+1}^{(j+1)}-\sigma^{(j+1)}_{j+1}\right)}}{\sum_{\ell=1}^j\gamma_\ell+\lambda
\mean \e{-\Dcal_{j}^k(\gamma)\left(\rho_{j}^{(j)}-\sigma^{(j)}_{j}\right)}}.
\]
The proof is finished after invoking (\ref{eq:transformlengths})
and noting that for $j=1,\ldots,k-1$,
\[
\frac{c_k}{c_j}\left[\lambda+\sum_{\ell=1}^j\gamma_\ell +\Dcal_j^k(\gamma)\right]=
\frac{c_k}{c_{j+1}}\left[\lambda+\sum_{\ell=1}^j\gamma_\ell +\Dcal_{j+1}^k(\gamma)\right]=
\lambda + \sum_{\ell=1}^{k-1} \gamma_\ell-\sum_{\ell=j+1}^{k-1}
\frac{c_k}{c_\ell}\gamma_\ell,
\]
and
\[
\frac{c_k}{c_1}\left[\lambda+\Dcal_1^k(\gamma)\right] =
\lambda + \sum_{\ell=1}^{k-1} \gamma_\ell-\sum_{\ell=1}^{k-1}
\frac{c_k}{c_\ell}\gamma_\ell,
\]
as the reader readily verifies.
\endproof

\subsection{A priority fluid system.}
\label{sec:priority}
In this subsection, we analyze a single station
which is drained at a constant rate $\mathtt{r}>0$.
It is fed by $n$ external inputs (`traffic classes') $J_1(t),\ldots,J_n(t)$,
each equipped with its own (infinite-capacity) buffer.
The queue discipline is (preemptive resume) {\it priority}, meaning that
for each $i=1,\ldots,n$,
the $i$-th buffer is continuously drained only if first $i-1$ buffers
do not require the full capacity $\mathtt{r}$.
We call such a system a {\it priority fluid system}.

The aim of this section is to find the Laplace transform of $(W,E)$,
where $W_j=W_j(\infty)$ is the stationary buffer content of class-$j$
input traffic, and $E_j=E_j(\infty)$ is the stationary age of the
busy period for class $j$.
We impose the following assumptions.
\begin{description}
\item[{\bf P1}] $J$ is an $n$-dimensional L\'evy
process with mutually independent components, and its L\'evy
measure is supported on $\R_+^n$, $J(0)=0$,
\item[{\bf P2}] $J_j(t)$ are nondecreasing for $j=2,\ldots,n$, and
\item[{\bf P3}] $J$ is integrable and $\sum_{i=1}^n \mean J_i(1)<\mathtt{r}$.
\end{description}

The central idea is that $W$ evolves in the same manner as
the solution to the Skorokhod problem that corresponds to
a tandem fluid network
$(J,r,P)$, with $r=(\mathtt{r},\ldots,\mathtt{r})'$ and $P=(p_{ij})$ such
that $p_{i,i+1}=1$ for $i=1,\ldots,n-1$ and $p_{ij}=0$ otherwise.
This equivalence has been noticed, for instance, by
Elwalid and Mitra~\cite{elwalidmitra:multiplexing1995}.
It allows us to use the notation of Section~\ref{sec:tandemgeneralities}.

It is important to observe that {\bf P1--P3} for the priority system
implies {\bf T1--T6} for the corresponding tandem fluid network,
except that the second inequality in {\bf T1} only holds as a weak inequality.
However, as remarked in Section~\ref{sec:tandemgeneralities},
the Laplace transform of the distribution of $(W,\widetilde B)$
is then still given in Theorem~\ref{th:tandem}.

The steady-state ages of the busy periods $E$ can also be expressed in terms
of the solution $(W,L)$ to this Skorokhod problem, but it does
{\em not always} equal $\widetilde B$ as in Section~\ref{sec:tandemgeneralities}.
To see this, notice that if class-$1$ traffic (highest priority) arrives
to an empty system at time $t$, we have $W_2(t)=0$, while $\widetilde W_2(t)>0$
so that $\widetilde B_2(t)>0$. However, it must hold that $E_2(t)=0$.

Still, the following theorem shows that it is possible to express the
distribution of $(W,E)$ in terms of $(W,\widetilde B)$.

\begin{Th}\label{th:priority}
Consider a priority fluid network for which {\bf P1--P3}
holds. Then for $\omega,\beta\in \R_+^n$, the transform
$\mean\e{-\langle\omega,W\rangle-\langle \beta,E\rangle}$ equals
\[
\mean \e{-\langle\omega,W\rangle-\langle \beta,\widetilde B\rangle}
+\sum_{j=2}^n \mean \left[\e{-\sum_{\ell=1}^{j-1} \omega_\ell W_\ell -
\sum_{\ell=1}^{j-1} \beta_\ell \widetilde B_\ell}
\left(1-\e{-\beta_j \widetilde B_j}\right);W_j=\ldots=W_n=0\right].
\]
\end{Th}
\proof{Proof}
In principle, $E_j$ equals $\widetilde B_j$, except when $W_j=0$.
In fact, it follows from the above reasoning that
\begin{eqnarray*}
\mean \e{-\langle \omega, W\rangle-\langle \beta, E\rangle} &=&
\mean \left[\e{-\omega_1 W_1-\beta_1 \widetilde B_1}; W_2=\ldots=W_n=0\right]\\
&& \mbox{}+ \sum_{j=2}^n
\mean \left[\e{-\sum_{\ell=1}^j\omega_\ell W_\ell-
\sum_{\ell=1}^j \beta_\ell \widetilde B_\ell}; W_j>0,
W_{j+1}=\ldots=W_n=0\right].
\end{eqnarray*}
Now use the fact that $\{W_j>0\}$ is the complement of $\{W_j=0\}$
and rearrange terms.
\endproof

\medskip
If the $J_2,\ldots,J_n$ are {\em strictly} increasing,
it can be seen (for instance with Theorem~\ref{th:tandem}) that
\[
\mean \left[\e{-\sum_{\ell=1}^{j-1} \omega_\ell W_\ell -
\sum_{\ell=1}^{j-1} \beta_\ell \widetilde B_\ell}
\left(1-\e{-\beta_j \widetilde B_j}\right);W_j=\ldots=W_n=0\right]=0.
\]
Therefore, in that case, we have the equality in distribution
$(W,E) =_{\rm d} (W,\widetilde B)$.

Another important special case is when $J_1,\ldots,J_n$ are compound
Poisson processes, say with intensities $\lambda_1,\ldots,\lambda_n$ respectively.
Much is known about the resulting priority system, see for instance
Jaiswal~\cite{jaiswal:priority1968} for this and related models.
To our knowledge, the distribution of $(W,E)$ has not been
investigated.
However, it is given by Theorem~\ref{th:priority} and Theorem~\ref{th:tandem}
upon noting that
$\theta^J_\ell(\omega)\to\lambda_\ell$ as $\omega\to\infty$.
Since it is not so instructive to write out the resulting formulas, we
leave this to the reader.

\section*{Acknowledgments}
The authors are grateful to Michel Mandjes, whose question prompted us to write this paper,
for his interest in this project.
KD, ABD and TR are partially supported by
KBN Grant No 1 P03A 031 28 (2005--2007), NWO Grant 631.000.002, and
KBN Grant No 2 P03A 020 23 (2002--2004) respectively.
Part of this work was done while ABD was with University College Cork, Ireland.

\appendix
\section{Appendix: some calculations for a compound Poisson process with negative drift}
\label{app:prelimcpdPsdrift}
In this appendix, we study a compound Poisson process $Z$
with negative drift, and
derive some results on the excursions of $Z-\underline Z$ from 0, just before
its entrance to 0.
These results are applied in Section~\ref{sec:distrH}.

Let us first fix the notation.
Throughout this appendix, $Z$ is a L\'evy process on $(\Omega, \mathcal F, \pr)$
with Laplace exponent
\[
\psi_{-Z}(\beta):=
\log \mean \e{-\beta Z(1)} = c \beta -\lambda \int_{\R_+}
\left(1-\e{-\beta z}\right) F(dz),
\]
where $c>0$, $\lambda\in(0,\infty)$, and $F$ is a probability distribution on $(0,\infty)$.
That is, $Z$ is a compound Poisson process under $\pr$ with rate $\lambda$ and
negative drift $-c$, and its (positive) jumps are governed by $F$.
We suppose that $\mean Z(1)<0$, so that $Z$ drifts to $-\infty$.
In analogy to Section~\ref{sec:distrH},
the inverse of $\psi_{-Z}$ is denoted by $\Phi_{-Z}$;
it is uniquely defined since $\psi_{-Z}$ is increasing. Observe that $\Phi_{-Z}(0)=0$.

Set $T_0=0$, and let $T_i$ denote the epoch of the $i$-th jump of $Z$.
To the $i$-th jump of $Z$,
we associate a vector of {\em marks}, denoted by $M_i\in\R^m_+$ (for some $m\in\Z_+$).
We suppose that $M_i$ is independent of the process $T\equiv\{T_n:n\geq 1\}$, and
that it is also independent of $(Z(T_j)-Z(T_j-), M_j)$ for $j\neq i$.
However, we allow for a dependency between $M_i$ and $Z(T_i)-Z(T_i-)$.
In fact, an interesting choice for $M_i$ is $M_i=Z(T_i)-Z(T_i-)$
(so that $m=1$).

Define $\tau_-$ as the first hitting time of zero, and
$N_-$ as the index of the last jump before $\tau_-$, i.e.,
\[
\tau_-:=\inf\{t\geq 0: Z(t)=0\},\quad
N_-=\inf\{n\geq 0: Z({T_{n+1}-})\leq 0\}.
\]
Write $\pr_\xi$ for the law of $Z+\xi$ under $\pr$ with initial mark $M_0=M$.
We suppose that the initial condition
$(\xi,M)$ is independent of $Z$, and has the same distribution as
$(Z(T_1)-Z(T_1-), M_1)$.
Observe that both $\tau_-$ and ${N_-}$ are $\pr_\xi$-almost surely finite,
and that (by the
Markov property) the `overshoot of the first excursion'
$T_{{N_-}+1}-\tau_-$ has an exponential distribution with parameter $\lambda$.

In this appendix, it is our aim to characterize the $\pr_\xi$-distribution
of $\tau_-$ (excursion length),
$\tau_--T_{N_-}$ (excursion `undershoot'), and $M_{N-}$ (mark of the last jump).
Overshoots and undershoots have been studied extensively in the literature.
However, as opposed to what we have here,
these results are all related to the situation that a
L\'evy process can cross a boundary by jumping over it
(strictly speaking, this is the only case where the terms `overshoot' and `undershoot'
seem to be appropriate). See Doney and
Kyprianou~\cite{doneykyprianou:overshootsundershoots2006} for a recent
contribution and for references.

In view of the
results of Dufresne and Gerber~\cite{dufresnegerber:surpluses1988}, 
it is tempting to believe that $\tau_--T_{N_-}$ has an exponential distribution.
However, it turns out that this `undershoot' has a completely different distribution.

\begin{Prop}
\label{prop:excdistr}
We have for $\beta,\gamma\geq 0$ and $\kappa\in\R^m_+$,
\begin{eqnarray*}
\mean_\xi \e{-\beta (\tau_--T_{N_-})-\gamma \tau_-
-\langle \kappa,M_{N_-}\rangle}&=&
\frac{\left[\beta + \gamma-
c\Phi_{-Z}(\gamma)+\lambda\right]\mean \e{-(\beta +\gamma + \lambda)\xi/c
-\langle \kappa,M\rangle}}
{\beta + \lambda \mean \e{-(\beta +\gamma + \lambda)\xi/c}} \\&=&
\frac{\left[\beta + \lambda\mean_\xi \e{-\gamma \tau_-}\right]
\mean \e{-(\beta +\gamma + \lambda)\xi/c -\langle \kappa,M\rangle}}
{\beta + \lambda \mean \e{-(\beta +\gamma + \lambda)\xi/c}}.
\end{eqnarray*}
\end{Prop}

To prove this proposition,
we need an auxiliary result on Poisson processes.
Consider a Poisson point process $N(t)$
with parameter $\mu$, and let
$\zeta$ be a positive random variable, independent of $N$.
Let $A(t)$ be the backward recurrence time process defined by $N$,
that is the time from $\zeta$ to the nearest point to the left.
The following lemma
characterizes the joint distribution of $N(\zeta)$, $A(\zeta)$, and $\zeta$.

\begin{Lemma}
\label{lem:overshootundershoot}
We have for $\beta,\gamma\geq 0$ and $0\leq s\leq 1$,
\[
\mean s^{N(\zeta)} \e{-\beta A(\zeta)-\gamma\zeta} =
\frac{\beta}{\beta+s\mu}\mean \e{-(\beta+\gamma+\mu)\zeta}
+\frac{s\mu}{\beta+s\mu}\mean \e{-[\gamma+(1-s)\mu]\zeta}.
\]
\end{Lemma}
\proof{Proof}
We only prove the claim for $\gamma=0$; the general case follows by
replacing the distribution of $\zeta$ by the (defective) distribution
of $\tilde \zeta$ given by
$\mean \e{-\beta \tilde \zeta}=\mean \e{-(\beta+\gamma)\zeta}$.
Let $U_0=0$ and $U_1,U_2,\ldots$ be the location of consecutive points of $N$.
Observe that
\begin{eqnarray}
\mean s^{N(\zeta)} \e{-\beta A(\zeta)} &=&
\sum_{n=0}^\infty s^n \mean \left[\e{-\beta (\zeta-U_n)}; 0\leq
\zeta-U_n\leq U_{n+1}-U_n\right]\nonumber\\
&=&\sum_{n=0}^\infty s^n \int_0^\infty \int_0^t \e{-(\beta+\mu)(t-x)} \pr_{U_n}(dx)\pr_\zeta(dt)
=\sum_{n=0}^\infty s^n \phi_n(\mu+\beta),\label{eq:sumphi}
\end{eqnarray}
where
\[
\phi_n(\beta):=\mean \left[\e{-\beta \left(\zeta-U_n\right)};
\zeta\geq U_n\right].
\]
Clearly, $\phi_0(\beta)=\mean \e{-\beta\zeta}$.
If we let $B$ be the forward recurrence time process, we have for $n\geq 1$,
\begin{eqnarray*}
\phi_n(\beta)&=&\mean\left[\e{-\beta(\zeta-U_n)}; \zeta\ge U_{n-1}\right]-
\mean \left[\e{-\beta (\zeta-U_n)};U_{n-1}\leq \zeta< U_n\right]\\
&=&\mean\left[\e{-\beta (\zeta-U_{n-1})+\beta(U_n-U_{n-1})};
\zeta\geq U_{n-1}\right]-\mean\left[\e{-\beta \left(\zeta-U_{n}\right)};
U_{n-1}\leq \zeta< U_n\right]\\
&=&\mean\left[\e{\beta(U_n-U_{n-1})}\right]\mean\left[\e{-\beta (\zeta-U_{n-1})};
\zeta\geq U_{n-1}\right]\\&&\mbox{}-\mean\left[\left.\e{\beta B(\zeta)}\right|
N(\zeta)=n-1\right]\pr(N(\zeta)=n-1)\\
&=& \frac{\mu}{\mu-\beta}\left[\phi_{n-1}(\beta)-\pr(N(\zeta)=n-1)\right],
\end{eqnarray*}
where we used the lack-of-memory property of the exponential
distribution for the last equality.
After iteration, we obtain
$$\phi_n(\beta)=\left(\frac{\mu}{\mu-\beta}\right)^n\mean \e{-\beta\zeta}
-\sum_{i=0}^{n-1}\left(\frac{\mu}{\mu-\beta}\right)^{n-i}\pr(N(\zeta)=i).$$
Therefore, taking $0<s<\beta/\mu$ (later we may use an analytic-continuation
argument), we deduce from (\ref{eq:sumphi}) that
\[
\mean\left[s^{N(\zeta)} \e{-\beta A(\zeta)}\right]=
\mean \e{-(\beta+\mu)\zeta}\sum_{n=0}^\infty
\left(-\frac{s\mu}{\beta}\right)^n-\sum_{n=1}^\infty s^n\sum_{i=0}^{n-1}
\left(-\frac{\mu}{\beta}\right)^{n-i}\pr(N(\zeta)=i).
\]
The double sum in this expression can be rewritten as
\[
-\frac{s\mu}{\beta+s\mu}\sum_{i=0}^\infty s^i \pr(N(\zeta)=i)\\
=-\frac{s\mu}{\beta+ s\mu}\mean \e{-(1-s)\mu\zeta},
\]
and the claim follows.
\endproof

\medskip
Lemma~\ref{lem:overshootundershoot} is the main ingredient to prove
Proposition~\ref{prop:excdistr}.
\medskip

\proof{Proof of Proposition~\ref{prop:excdistr}}
The crucial yet simple observation is that
\begin{eqnarray}
\lefteqn{\mean_\xi \e{-\beta (\tau_--T_{N_-})-\gamma \tau_--\langle \kappa,M_{N_-}\rangle}}\nonumber\\
&=&\mean_\xi \left[\e{-\beta (\tau_--T_{N_-})-\gamma \tau_--\langle \kappa,M_{N_-}\rangle};
{N_-}=0\right]+
\mean_\xi \left[\e{-\beta (\tau_--T_{N_-})-\gamma \tau_--\langle \kappa,M_{N_-}\rangle};
{N_-}\geq 1\right] \nonumber \\
&=&
\mean \e{-(\lambda+\beta+\gamma)\xi/c -\langle \kappa,M\rangle}
+\mean_\xi \left[\e{-\beta (\tau_--T_{N_-})-\gamma \tau_-
-\langle \kappa,M_{N_-}\rangle}; {N_-}\geq 1 \right].
\label{eq:renewaltype}
\end{eqnarray}
To analyze the second term, we exploit the fact that there are several
excursions of $Z-\underline Z$ from 0.
Therefore, we set
\[
C(t):=\inf\{s\geq 0: Z(s)-Z(0)=-t\},
\]
where an infimum over an empty set should be interpreted as infinity.

It is obvious that $C$ is a subordinator with drift $1/c$, and that it
jumps at rate $\lambda/c$ with jumps distributed as
$\tau_-$ under $\pr_\xi$. This observation implies with
Theorem~VII.1 of Bertoin~\cite{bertoin:levy1996} that
\begin{equation}
\label{eq:characPhi}
\Phi_{-Z}(\gamma)=\frac\gamma c+\frac\lambda c \left(1-\mean_\xi \e{-\gamma \tau_-}\right).
\end{equation}

Lemma~\ref{lem:overshootundershoot} can be applied to
the Poisson process $N$ constituted by the jump epochs of $C$, $\mu=\lambda/c$, and
$\zeta=\xi$.
Each jump of $C$ corresponds to an excursion of $Z-\underline Z$ from 0,
for which the `excursion overshoot', the excursion length, and the marks
of the last jump are of interest.
Observe that these quantities have the same distribution
as $\tau_--T_{N_-}$, $\tau_-$, and $M_{N_-}$ respectively.
Using the notation of Lemma~\ref{lem:overshootundershoot}, this yields
\begin{eqnarray}
\lefteqn{\mean_\xi \left[\e{-\beta (\tau_--T_{N_-}) -\gamma\tau_--\langle \kappa,M_{N_-}\rangle};
{N_-}\geq 1 \right]} \nonumber\\&=&
\mean \left[\left(\mean_\xi \e{-\gamma\tau_-}\right)^{N(\xi)-1}
\e{-\beta A(\xi)/c-\gamma\xi/c}; N(\xi)\geq 1\right]\mean_\xi
\e{-\beta (\tau_--T_{N_-})-\gamma \tau_--\langle \kappa,M_{N_-}\rangle}.\label{eqn:identity}
\end{eqnarray}
Therefore, Lemma~\ref{lem:overshootundershoot} yields
\begin{eqnarray*}
\mean \left[s^{N(\xi)-1} \e{-\beta A(\xi)/c-\gamma\xi/c}; N(\xi)\geq 1\right]&=&
\frac{\mean \left[s^{N(\xi)} \e{-\beta A(\xi)/c-\gamma\xi/c}\right]-
\mean \e{-(\lambda+\beta+\gamma)\xi/c}}{s}\\
&=& \frac\lambda{\lambda s+\beta} \left[\mean \e{-((1-s)\lambda+\gamma)\xi/c}-
\mean \e{-(\lambda+\beta+\gamma)\xi/c}\right].
\end{eqnarray*}
Upon combining this with (\ref{eq:renewaltype}) and (\ref{eqn:identity}), we arrive at
\[
\mean_\xi \e{-\beta (\tau_--T_{N_-})-\gamma \tau_--\langle \kappa,M_{N_-}\rangle} = \frac{
\left[\beta+\lambda\mean_\xi \e{-\gamma\tau_-}\right]
\mean \e{-(\lambda+\beta+\gamma)\xi/c-\langle \kappa,M\rangle}}
{\lambda\mean_\xi \e{-\gamma\tau_-}+\beta -\lambda
\mean \e{-(\lambda(1-\mean_\xi \e{-\gamma\tau_-})+\gamma)\xi/c}+
\lambda \mean \e{-(\lambda+\beta+\gamma)\xi/c}},
\]
which, with the help of (\ref{eq:characPhi}), reduces to
\[
\frac{\left[\beta + \gamma-
c\Phi_{-Z}(\gamma)+\lambda\right]\mean \e{-(\beta +\gamma + \lambda)\xi/c-
\langle \kappa,M\rangle}}
{\beta + \gamma-c\Phi_{-Z}(\gamma)-\lambda\left(
\mean \e{-\Phi_{-Z}(\gamma)\xi}-1\right)+\lambda
\mean \e{-(\beta +\gamma + \lambda)\xi/c}}.
\]
By definition of $\Phi_{-Z}$, we have
\[
\gamma = \psi_{-Z}(\Phi_{-Z}(\gamma)) = c\Phi_{-Z}(\gamma)+\lambda\left(
\mean \e{-\Phi_{-Z}(\gamma)\xi}-1\right),
\]
and the claim follows.
\endproof

\footnotesize
\bibliography{y:/bibdb}
\bibliographystyle{amsplain}

\end{document}